\RequirePackage{amsthm} 
\documentclass[pdflatex, iicol, sn-mathphys-num]{sn-jnl}
\usepackage{graphicx}%
\usepackage{multirow}%
\usepackage{amsmath,amssymb,amsfonts}%
\usepackage{amsthm}%
\usepackage{mathrsfs}%
\usepackage[title]{appendix}%
\usepackage{xcolor}%
\usepackage{textcomp}%
\usepackage{manyfoot}%
\usepackage{booktabs}%
\usepackage{algorithm}%
\usepackage{algorithmicx}%
\usepackage{algpseudocode}%
\usepackage{listings}%
\usepackage{bm}
\usepackage{enumitem}
\usepackage{cleveref}
\newcommand{\argmin}{\mathop{\mathrm{argmin}}}

\DeclareMathOperator{\dd}{d\!}

\usepackage{tcolorbox}

\usepackage{subfig}
\usepackage{tikzscale}
\usepackage{pgfplots}
\usepackage{tabularray}
\pgfplotsset{width=7cm,compat=1.8}
\usetikzlibrary{pgfplots.groupplots}{}

\definecolor{color0}{rgb}{0.7843, 0.7843, 0.7843}
\definecolor{color1}{rgb}{0, 0.4470, 0.7410}
\definecolor{color2}{rgb}{0.8500, 0.3250, 0.0980} \definecolor{color3}{rgb}{0.9290, 0.6940, 0.1250}
\definecolor{color4}{rgb}{0.7060, 0.3840, 0.7650}
\definecolor{color5}{rgb}{0.4660, 0.6740, 0.1880}
\definecolor{color6}{rgb}{0.3010, 0.7450, 0.9330}
\definecolor{color7}{rgb}{0.6350, 0.0780, 0.1840}
\definecolor{color8}{rgb}{0.0, 0.4078, 0.3412}

\allowdisplaybreaks

\theoremstyle{thmstyleone}%
\newtheorem{theorem}{Theorem}
\newtheorem{proposition}[theorem]{Proposition}

\theoremstyle{thmstyletwo}%
\newtheorem{remark}{Remark}%

\theoremstyle{thmstylethree}%

\begin{document}

\title[SiMPL]{A Simple Introduction to the SiMPL Method for Density-Based Topology Optimization}

\author[1]{\fnm{Dohyun} \sur{Kim}}\email{dohyun\_kim@brown.edu}
\author[2]{\fnm{Boyan S.} \sur{Lazarov}}\email{lazarov2@llnl.gov}
\author[3]{\fnm{Thomas M.} \sur{Surowiec}}\email{thomasms@simula.no}
\author[1]{\fnm{Brendan} \sur{Keith}}\email{brendan\_keith@brown.edu}
\affil[1]{\orgdiv{Division of Applied Mathematics}, \orgname{Brown University}, \orgaddress{\city{Providence}, \postcode{02912}, \state{RI}, \country{United States of America}}}
\affil[2]{\orgname{Lawrence Livermore National Laboratory}, \orgaddress{\city{Livermore}, \postcode{94550}, \state{CA}, \country{United States of America}}}
\affil[3]{\orgdiv{Department of Numerical Analysis and Scientific Computing}, \orgname{Simula Research Laboratory}, \orgaddress{\city{Oslo}, \postcode{0164}, \country{Norway}}}

\abstract{
  We introduce a novel method for solving density-based topology optimization problems: \underline{Si}gmoidal \underline{M}irror descent with a \underline{P}rojected \underline{L}atent variable (SiMPL). The SiMPL method (pronounced as ``the simple method'') optimizes a design using only first-order derivative information of the objective function. The bound constraints on the density field are enforced with the help of the (negative) Fermi--Dirac entropy, which is also used to define a non-symmetric distance function called a Bregman divergence on the set of admissible designs. This Bregman divergence leads to a simple update rule that is further simplified with the help of a so-called latent variable. 
  Because the SiMPL method involves discretizing the latent variable, it produces a sequence of pointwise-feasible iterates, even when high-order finite elements are used in the discretization. Numerical experiments demonstrate that the method outperforms other popular first-order optimization algorithms. To outline the general applicability of the technique, we include examples with (self-load) compliance minimization and compliant mechanism optimization problems.
}


\maketitle

\section{Introduction}
Topology optimization (TO) is a design process that seeks to find the optimal distribution of material in a selected physical domain subject to multiple constraints. The method has been utilized in many engineering applications, varying in scale and complexity, e.g., microstructures, such as the design of photonic crystals, and exotic meta-materials or macro-structures, such as aircraft wing designs, bridges, buildings, and mechanical assemblies \cite{Bendsoe2004}. Two main approaches can be distinguished in the literature for representing the design: the level-set approach and the density-based approach. As the name suggests, the level-set approach \cite{allaire2004level} relies on a level-set function $\phi$ to represent the material interface $\left(\phi=0\right)$. On the other hand,  the material distribution in density-based approaches (the main focus of this paper) utilizes a density field $0\leq\rho\leq 1$, taking the value zero in the void regions and one in the subdomain occupied with solid material.

\paragraph*{Density variables}
Gradient-based optimization algorithms for density-based TO require continuous and smooth transitions between void and solid regions. Thus, the actual physical density is represented via a series of transformations applied to the original density field $\rho$. The most popular formulation utilizes the original density field, a filtered field, and, finally, the physical density field \cite{Lazarov2016}. The filtered field $\tilde\rho$ is obtained by convolving the original density field with a filter function. This step regularizes the optimization problem and guarantees the existence of a solution  \cite{bourdin2001}. Subsequently, the filtered field is mapped to the physical density field through a material interpolation model such as the solid isotropic material with penalization (SIMP) \cite{Bendsoe1999} or rational approximation of material property (RAMP) \cite{stolpe2001} models.  

\paragraph*{Optimization methods}
Given a suitable spatial discretization, the original topology optimization problem becomes a finite-dimensional constrained optimization problem that can be solved with gradient-based, black box optimization solvers. For density-based TO, the most popular optimization methods are the Optimality Criteria method (OC) \cite{Bendsoe2004} and the Method of Moving Asymptotes (MMA) \cite{svanberg1987mma}. The simplicity and efficiency of OC, demonstrated in publicly available MATLAB  code provided in \cite{sigmund2001matlab, Andreassen2011}, have made it a popular choice for TO problems. However, the update rule in OC is heuristic, and the theoretical foundation of the method is still under development \cite{ananiev2005ocmpgd}. The main alternative, MMA, solves a sequence of subproblems to find a stationary point of the optimization problem. The method can handle more general constraints than OC, and a globally convergent version is available \cite{zillober1993gcmma}.
However, MMA is a general-purpose finite-dimensional optimization algorithm involving many complex parameters that are difficult for an average user to tune.
Moreover, many implementations show mesh-dependent behavior, i.e., 
the number of optimization steps increases after mesh refinement when starting from the same initial guess.
The same problem appears with the OC method, or any other algorithm implemented without the correct derivative-to-gradient transformation (i.e., Riesz map) \cite{Schwedes2017, Petra2023}.
For an alternative second method, we refer a trust-region based method \cite{Kouri2021}.

\paragraph*{Optimize-then-discretize}
Mesh-dependence in optimization can often be avoided by taking an optimize-then-discretize approach, which is utilized in \Cref{sub:opt-disc} of the present article.
This involves first deriving a theoretical optimization algorithm for the non-discretized problem and then discretizing that algorithm to obtain the enacted (i.e., practical) version. This approach is often used in level-set topology and shape optimization \cite{allaire2004level,Bergmann2019} and examples in density-based topology optimization can be found in \cite{Bendsoe2004, Evgrafov2014,PapoutsisKiachagias2016,Jensen2017}. The optimize-then-discretize paradigm can make it easier to exploit the structure of the original infinite-dimensional optimization problem. Generally, it also provides greater freedom for selecting discretizations and opens the possibly for adapting the discretization between algorithm iterations. 

\paragraph*{Preserving bound constraints}
In this work, we consider an approach inspired by the proximal Galerkin method for variational inequalities introduced in \cite{keith2023proximal}. The proposed method is a so-called mirror descent \cite{nemirovskij1983problem,beck2003mirror,teboulle2018nolip} algorithm tailored to the precise mathematical structure present in density-based TO. In turn, it produces a sequence of feasible design iterates regardless of the polynomial order of the underlying finite element discretization, resulting in an optimized density field guaranteed to satisfy the bound constraints at every point in the design domain.
In particular, the SiMPL method utilizes a latent variable representation of the orginal density function $\rho$.
This simplifies the update rule and ensures that the discretized density field is always updated in a bound-preserving manner.

\paragraph*{Innovations}
Initial numerical experiments in \cite{keith2023proximal} revealed sensitivities to the choice of step size, affecting both the convergence and quality of the optimized solution.
We overcome this issue through an adaptive step size strategy closely related to the Barzilai--Borwein method \cite{barzilai1988bbm} and the techniques discussed in \cite{bauschke2017-NoLip,  teboulle2018nolip, lu2019relative}. The estimated step size is used as an initial guess for a line search algorithm that ensures a monotonically decreasing objective function value. In our numerical experiments, the number of iterations is reduced significantly compared to the linearly growing step size rule utilized in \cite{keith2023proximal}.

\paragraph*{Software availability}
An open-source implementation of the SiMPL method is available in the finite element library MFEM \cite{tzanio2010mfem,andrej2024mfem} (Example 37).
This implementation accompanies \cite{mfem_simpl}, which contains all of the code to reproduce our numerical experiments.

\paragraph*{Outline}
The rest of the paper is organized as follows. We begin by using the discretize-then-optimize paradigm to derive the SiMPL method for a minimum compliance problem discretized with lowest-order finite elements.
This facilitates an easier transition to the infinite-dimensional formulation derived in \Cref{sub:opt-disc}.
We then describe a backtracking line search algorithm, together with two choices of sufficient decrease conditions, that we recommend to be used with the SiMPL method.
After the full algorithm derivation, we proceed with numerical experiments on a 2D MBB beam problem to verify mesh-independent convergence of the SiMPL method and provide performance comparisons to OC and MMA. Finally, to demonstrate the general applicability of the SiMPL method, we showcase optimized solutions to self-weight compliance minimization and compliant mechanism problems.
\section{The SiMPL method}\label{sec:mirror-descent}
The SiMPL method is based on mirror descent \cite{nemirovskij1983problem}, which is a non-Euclidean generalization of the well-known steepest descent method leveraging a type of squared distance function called a Bregman divergence \cite{Bregman1967}.
The specific form of mirror descent we are proposing is tailored to the bound constraints $0\leq \rho\leq 1$; cf.\ \eqref{eq:bound-constraint}, below.
To streamline the initial exposition, we begin by using the discretize-then-optimize paradigm to derive the SiMPL method in the most common setting for density-based topology optimization: piecewise-constant (lowest-order) discrete densities on grid-like meshes.
\Cref{sub:LineSearch} then introduces various step size selection strategies for optimized efficiency, while \Cref{sub:stopping_criteria} proposes some convenient stopping criteria for the SiMPL method.
Finally, using an optimize-then-discretize approach supported by mathematical results from \cite{simplmath}, \Cref{sub:opt-disc} describes how the SiMPL method can be applied to high-order discretizations on more generally-meshed domains.

\subsection{First discretize then optimize}
\label{sub:FDTO}

We now derive the SiMPL method for piecewise-constant discrete densities on grid-like meshes using the discretize-then-optimize paradigm.

\paragraph{Problem definition}
We focus on topology optimization of a linearly elastic structure in a design domain $\Omega \subset \mathbb{R}^d$ partitioned into Cartesian cells $\Omega_h=\{\Omega_1,...,\Omega_{N_\rho}\}$ with maximum diameter $h > 0$.
We seek finite-dimensional approximations of the density ${\rho} \in L^2(\Omega)$, filtered density ${\tilde{\rho}} \in H^1(\Omega)$, and displacement ${u} \in V\subset [H^1(\Omega)]^d$ in $Q_h$, $\tilde{Q}_h$, and $V_h$, respectively, where
\begin{subequations}
\label{eq:LowestOrderSpaces}
  \begin{align}
    Q_h         &:=\{q\in L^2(\Omega):q|_{\Omega_i}\in Q_0(\Omega_i)\;\forall \Omega_i\in\Omega_h\},                 \\
    \tilde{Q}_h &:=\{\tilde{q}\in H^1(\Omega):\tilde{q}|_{\Omega_i}\in Q_1(\Omega_i)\;\forall \Omega_i\in\Omega_h\}, \\
    V_h         &:=\{v\in V:v|_{\Omega_i}\in [Q_1(\Omega_i)]^d\; \forall \Omega_i\in\Omega_h\},
  \end{align}
\end{subequations}
and each $Q_k(\Omega_i)$, with $k = 0$ or $1$, is the space of constant functions or multilinear polynomials, respectively, defined over the cell $\Omega_i$.
We can represent functions in the finite-dimensional spaces~\eqref{eq:LowestOrderSpaces} by coefficient vectors containing the functions' cell/nodal values.
These coefficient vectors are denoted with bold symbols by $\bm{\rho}\in\mathbb{R}^{N_\rho}$, $\tilde{\bm{\rho}}\in \mathbb{R}^{N_{\tilde{{\rho}}}}$, $\mathbf{u}\in\mathbb{R}^{N_u}$, where $N_\rho$, $N_{\tilde{{\rho}}}$, and $N_u$ are the numbers of the degrees of freedom of the discretized density, filtered density, and displacement field, respectively.
Following \cite{Bendsoe2004, lazarov2011-filter, Andreassen2011}, the discretized topology optimization problem can now be written as follows:~
\begin{subequations}
  \label{eq:disc-prob}
  \begin{align}
    \label{eq:objective}
    \min_{\bm{\rho}, \tilde{\bm{\rho}}, \mathbf{u}}\  & \widehat{F}(\tilde{\bm{\rho}}, \mathbf{u})                                  \\
    \text{subject to }
    \label{eq:state-eq}
                                                      & \mathbf{K}(\tilde{\bm{\rho}})\mathbf{u}=\mathbf{f},                         \\
    \label{eq:filter-eq}
                                                      & (\epsilon^2 \mathbf{A} + {\tilde{\mathbf{M}}})\tilde{\bm{\rho}}={\mathbf{N}}\bm{\rho}, \\
    \label{eq:bound-constraint}
                                                      & \mathbf{0}\leq \bm{\rho} \leq \mathbf{1},                                                                \\
    \label{eq:volume-constraint}
                                                      & \mathbf{1}^\top\mathbf{M}\bm{\rho}\leq\theta |\Omega|.
  \end{align}
\end{subequations}
Here, $\mathbf{K}$ is the linear elasticity tangent (stiffness) matrix, $\mathbf{A}$ is the stiffness matrix corresponding to the diffusion operator in a discretized PDE-filter \cite{lazarov2011-filter}, ${\tilde{\mathbf{M}}}$ is the (symmetric) mass matrix for the filtered density variable, ${\mathbf{N}}$ is the (non-symmetric) mass matrix between the density and filtered density spaces, and $\mathbf{M}$ is the (diagonal) mass matrix for the density variable obeying $\mathbf{1}^\top\mathbf{M}\mathbf{1} = \operatorname{tr}\mathbf{M} = |\Omega|$.
$\widehat{F}(\tilde{\bm{\rho}}, \mathbf{u})$ is a prescribed objective function, $\mathbf{f}$ is a vector representation of a load applied to the system,  $\epsilon=r_{\min}/2\sqrt{3}$ is the filter coefficient with the minimum length scale $r_{\min} > 0$, and $0 < \theta < 1$ is a prescribed volume fraction.
Here and throughout, vector inequalities such as~\eqref{eq:bound-constraint} are understood component-wise.
In~\eqref{eq:state-eq}, the individual element contributions to the tangent matrix $\mathbf{K}$ are calculated as $\mathbf{K}_i=E_i\mathbf{K}_0$ where $\mathbf{K}_0$ is the element stiffness matrix for a unit stiffness, and $E_i$ is the material stiffness obtained by using the so-called solid isotropic material interpolation with penalization (SIMP) law \cite{Bendsoe1989} written as
\begin{equation*}
  E_i = E_{\min} + r(\tilde{\rho}_i)(E_{\max} - E_{\min})
  \,.
\end{equation*}
In this expression, $E_{\max}$ and $E_{\min}$ are the stiffnesses of the solid and void phases, respectively, and $r(\tilde{\rho}_i)=\tilde{\rho}_i^p$ is the physical density with the penalization exponent $p>1$.
Typically, the exponent is selected to be $p=3$ and $E_{\min}/E_{\max}=10^{-6}$.
In the following, we employ the reduced space approach by letting
\begin{equation}
\label{eq:ObjectiveFunction}
  F(\bm{\rho}) = \widehat{F}(\tilde{\bm{\rho}}(\bm{\rho}), \mathbf{u}(\tilde{\bm{\rho}}(\bm{\rho})))
  \,.
\end{equation}
We also use the symbol
\begin{equation}
\label{eq:AdmissibleSet}
  \mathcal{A}_h = \{ \bm{\rho} \mid \mathbf{1}^\top\mathbf{M}\bm{\rho}\leq\theta|\Omega|,~ \mathbf{0}\leq \bm{\rho}\leq \mathbf{1} \}
\end{equation}
to denote the set of admissible design vectors.
Together, this notation allows us to rewrite~\eqref{eq:disc-prob} as
\begin{equation}
\label{eq:ReducedProblem}
  \min_{\bm{\rho} \in \mathcal{A}_h}
  F(\bm{\rho})
  \,.
\end{equation}

\paragraph{Steepest descent}
To motivate the SiMPL method \eqref{eq:full-update}, we begin by recalling the steepest descent method under a weighted inner product.
Given the previous iteration $\bm{\rho}_k$ and a step size $\alpha_k > 0$, steepest descent finds its next iterate $\bm{\rho}_{k+1}$ by minimizing the following local quadratic approximation to $F(\bm{\rho})$:
\begin{multline}\label{eq:GD-obj}
  J(\bm{\rho};\bm{\rho}_k)=F(\bm{\rho}_k) +\mathbf{d}F_k^\top(\bm{\rho}-\bm{\rho}_k) \\+\frac{1}{2\alpha_k}(\bm{\rho}-\bm{\rho}_k)^\top\mathbf{B}(\bm{\rho}-\bm{\rho}_k).
\end{multline}
Here, $\mathbf{d}F_k\in\mathbb{R}^{N_\rho}$ denotes the $\ell^2$-gradient of $F(\bm{\rho}_k)$, satisfying $(\mathbf{d}F_k)_i=\partial (F(\bm{\rho}_k))/\partial (\bm{\rho}_k)_i$ for each $i = 1, 2, \ldots, N_{\rho}$, $\mathbf{B}\in\mathbb{R}^{N_\rho\times N_\rho}$ is symmetric positive-definite matrix, and $\alpha_k > 0$ is a prescribed step size.
The standard steepest descent method takes $\mathbf{B}=\mathbf{I}$, which in~\eqref{eq:GD-obj} corresponds to using the square of the $\ell^2$-norm to penalize $\bm{\rho}-\bm{\rho}_k$.
Notably, if we let $\mathbf{B}=\mathbf{M}$, this corresponds to steepest descent with a discrete $L^2(\Omega)$-norm.

\paragraph{Choice of inner product}
In the following, we use only $\mathbf{B}=\mathbf{M}$.
This choice makes our derivation consistent with the units and length scales of the underlying problem and the optimize-then-discretize derivation that is given in \Cref{sub:opt-disc}.
To this end, we denote the $L^2(\Omega)$-gradient vector,
\begin{align}
\label{eq:gradient}
  \mathbf{g}_k=\mathbf{M}^{-1}\mathbf{d}F_k,
\end{align}
which is directly related to the $L^2(\Omega)$-gradient function introduced later on, in~\eqref{eq:derivative}; see also~\eqref{eq:DiscreteGradientEquation}.

\paragraph{Steepest descent with projection}
Upon setting $\mathbf{B} = \mathbf{M}$ in~\eqref{eq:GD-obj}, a straightforward computation shows that
\begin{equation}
\label{eq:PGD}
  \begin{aligned}
    \bm{\rho}_{k+1}
     & =\argmin_{\bm{\rho}\in \mathcal{A}_h} J(\bm{\rho};\bm{\rho}_k)                                                 \\
     & =\min\{\mathbf{1},\max\{\mathbf{0},\bm{\rho}_k-\alpha_k\mathbf{g}_k-\alpha_k\mu_{k+1}\mathbf{1}\}\}
     \\
     & =\mathcal{P}(\bm{\rho_k} - \alpha_k\mathbf{g}_k)
     \,.
  \end{aligned}
\end{equation}
Here, $\mu_{k+1} \geq 0$ is a scalar Lagrange multiplier ensuring that the volume constraint~\eqref{eq:volume-constraint} is satisfied and $\mathcal{P}$ denotes the discrete $L^2(\Omega)$-projection onto the admissible set~\eqref{eq:AdmissibleSet}:
\begin{equation}
\label{eq:GD-proj}
  \mathcal{P}(\bm{\rho})
  =
  \argmin_{\mathbf{q} \in \mathcal{A}_h}
  \frac{1}{2}(\bm{\rho}-\mathbf{q})^\top\mathbf{M}(\bm{\rho}-\mathbf{q})
  \,.
\end{equation}
Using standard KKT-theory from nonlinear programming, see e.g., \cite[Chap. 12]{nocedal1999numerical}, and rearranging terms, we arrive at \eqref{eq:PGD}.
Mirror descent methods \cite{nemirovskij1983problem}, such as the SiMPL method, replace the weighted inner products in \eqref{eq:GD-obj} and \eqref{eq:GD-proj} with a Bregman divergence $D_\varphi(\cdot,\cdot)$.

\paragraph{Bregman divergences}
Introduced in 1967 by L.M.\ Bregman \cite{Bregman1967}, Bregman divergences generalize the concept of a squared Euclidean distance using the error in the linear approximation to a strictly convex function.
More specifically, the Bregman divergence induced by a strictly convex proper function $\varphi \colon \mathbb{R}^{N_\rho} \to \mathbb{R} \cup \{+\infty\}$ is given by:
\begin{equation}
\label{eq:BregmamDivergence}
  D_\varphi(\bm{\rho},\mathbf{q}) = \varphi(\bm{\rho}) - \varphi(\mathbf{q}) - \mathbf{d}\varphi(\mathbf{q})^\top(\bm{\rho} - \mathbf{q}),
\end{equation}
for all admissible $\bm{\rho},\bm{q} \in \mathbb{R}^{N_\rho}$.
If $\varphi$ is the weighted inner product $\varphi(\bm{\rho}) = \frac{1}{2}\bm{\rho}^\top \mathbf{B} \bm{\rho}$, then its Bregman divergence is exactly the squared distance function $\frac{1}{2}(\bm{\rho} - \mathbf{q})^\top \mathbf{B} (\bm{\rho} - \mathbf{q})$ appearing in \eqref{eq:GD-obj}.
Replacing this squared distance function with the general Bregman divergence~\eqref{eq:BregmamDivergence} and removing the remaining constant terms, we arrive at
\begin{equation}
  J_\varphi(\bm{\rho};\bm{\rho}_k) =
  \mathbf{d}F_k^\top \bm{\rho} + \frac{1}{\alpha_k}D_\varphi(\bm{\rho},\bm{\rho}_k)
  \,.
\end{equation}
Likewise, replacing the weighted inner product in~\eqref{eq:GD-proj} with a Bregman divergence allows us to define the so-called Bregman projection \cite{bauschke1997legendre}:
\begin{equation}
\label{eq:MD-proj}
  \mathcal{P}_{\varphi}(\bm{\rho})
  =
  \argmin_{\mathbf{q} \in \mathcal{A}_h}
  D_\varphi(\mathbf{q},\bm{\rho})
  \,.
\end{equation}

\begin{figure*}[t]
  \centering
  \includegraphics[width=0.75\textwidth]{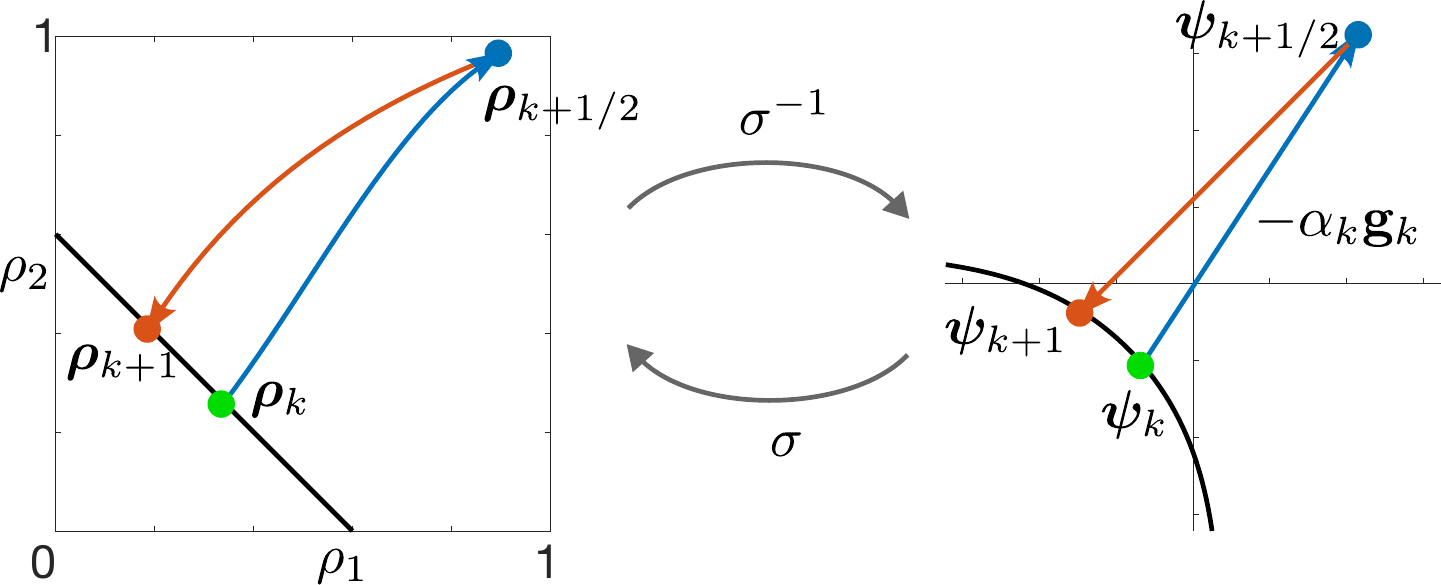}
  \caption{
    Schematics of the SiMPL method in primal (left) and latent (right) spaces, respectively, in $\mathbb{R}^2$.
    $\bm{\rho}_{k+1/2}=\sigma(\bm{\psi}_{k+1/2})=\sigma\big(\bm{\psi}_k-\alpha_k\nabla F(\bm{\rho}_k)\big)$ is an auxiliary step before the volume correction.
    Black curves represent the feasible set $K$ with volume constraint; cf.\ \eqref{eq:full-update}.
    Both the gradient step and volume correction are linear operations in the latent space (right), but are nonlinear in the primal space (left).
    \label{fig:mapping}
  }
\end{figure*}

\paragraph{Fermi--Dirac entropy}
To solve the topology optimization problem~\eqref{eq:disc-prob}, we are interested in a particular choice of $\varphi$, namely, the (negative) Fermi--Dirac entropy:
\begin{equation}\label{eq:disc-varphi}
  \varphi(\bm{\rho}) = \ln(\bm{\rho})^\top\mathbf{M}\bm{\rho}+\ln(\mathbf{1}-\bm{\rho})^\top\mathbf{M}(\mathbf{1}-\bm{\rho}).
\end{equation}
Here, the logarithms are understood to be applied component-wise.
We argue that this function encodes the geometry of the box constraint $\mathbf{0} \leq \bm{\rho} \leq \mathbf{1}$ appearing in~\eqref{eq:bound-constraint}.
In particular, its gradient is the inverse of the (logistic) sigmoid function, $\sigma(x)=1/(1+\exp(x))$, as can be readily verified:
\begin{equation}
\label{eq:DiscreteFermiDirac}
  \mathbf{M}^{-1}\mathbf{d}\varphi(\bm{\rho}) =\ln\left(\bm{\rho}/(\mathbf{1}-\bm{\rho})\right)= \sigma^{-1}(\bm{\rho}).
\end{equation}
This gradient is a structure-preserving mapping between the open set $(0,1)^{N_\rho}$ and $\mathbb{R}^{N_\rho}$; see \Cref{fig:mapping}.

\paragraph{Sigmoidal mirror descent with projection}
Assuming $\mathbf{0} < \bm{\rho}_k < \mathbf{1}$, the next mirror descent iterate can be obtained by minimizing $J_\varphi(\bm{\rho};\bm{\rho}_k)$ over $\bm{\rho}$ in the admissible set~\eqref{eq:AdmissibleSet}:
\begin{equation}
\label{eq:PMD}
  \begin{aligned}
    \bm{\rho}_{k+1}
     & =\argmin_{\bm{\rho}\in \mathcal{A}_h} J_\varphi(\bm{\rho};\bm{\rho}_k)                                              \\
     & =\sigma\big(\sigma^{-1}(\bm{\rho}_k)-\alpha_k\mathbf{g}_k-\alpha_k\mu_{k+1}\mathbf{1}\big)
     \\
     & =\mathcal{P}_{\varphi}\Big(\sigma\big(\sigma^{-1}(\bm{\rho}_k) - \alpha_k\mathbf{g}_k\big)\Big)
     .
  \end{aligned}
\end{equation}
Notice that, like~\eqref{eq:PGD}, $\mu_{k+1} \geq 0$ is a Lagrange multiplier corresponding to the volume constraint \eqref{eq:volume-constraint}.
However, unlike~\eqref{eq:PGD}, the middle expression does not involve clipping.
Instead, the Bregman projection $\mathcal{P}_{\varphi}$ in the third line of~\eqref{eq:PMD} is a smooth operator.

\paragraph{The latent variable}
If $\mathbf{0} < \bm{\rho}_0 < \mathbf{1}$, then by~\eqref{eq:PMD}, every subsequent iterate $\bm{\rho}_k$ satisfies the box constraint~\eqref{eq:bound-constraint} strictly (i.e., $\mathbf{0} < \bm{\rho}_k < \mathbf{1}$ for all $k \geq 1$).
Moreover, $\bm{\psi}_k := \sigma^{-1}(\bm{\rho}_k) \in \mathbb{R}^{N_\rho}$ is always a well-defined vector.
The SiMPL method chooses to evolve the latent variable $\bm{\psi}_k$ throughout the optimization process.
This is done for three main reasons:~
\medskip

  \noindent \textit{\textbf{First,}} the latent variable form of the update rule~\eqref{eq:PMD} has the following simple and convenient two-stage structure:
  \begin{subequations}\label{eq:full-update}
    \begin{align}
      \label{eq:gradient-update}
      \bm{\psi}_{k+1/2} & =\bm{\psi}_{k}-\alpha_k\mathbf{g}_k,                                  \\
      \label{eq:volume-update}
      \bm{\psi}_{k+1}   & =\bm{\psi}_{k+1/2}-\alpha_k\mu_{k+1}\mathbf{1}
      .
    \end{align}
  Stage one \eqref{eq:gradient-update} corresponds to an unconstrained gradient step in the latent space $\mathbb{R}^{N_\rho}$.
  Next, stage two \eqref{eq:volume-update} uniformly translates each component of the intermediary latent variable $\bm{\psi}_{k+1/2}$ until the volume constraint
  \[
    \mathbf{1}^\top\bm{M} \sigma(\bm{\psi}_{k+1/2}-\alpha_k\mu_{k+1}\mathbf{1})\leq\theta|\Omega|
  \]
  is satisfied.
  Note that we set $\mu_{k+1} = 0$ if $\mathbf{1}^\top\mathbf{M} \sigma(\bm{\psi}_{k+1/2}) < \theta|\Omega|$.
  Otherwise, we find the unique $\mu_{k+1} \geq 0$ solving the nonlinear equation
  \begin{equation}
  \label{eq:vol-correction}
    \mathbf{1}^\top\mathbf{M} \sigma(\bm{\psi}_{k+1/2}-\alpha_k\mu_{k+1}\mathbf{1})
    =
    \theta|\Omega|
    \,.
  \end{equation}
  \end{subequations}
  \Cref{eq:gradient-update,eq:volume-update} are linear update rules that are simple to implement.
  Meanwhile, solving \eqref{eq:vol-correction} for $\mu_{k+1}$ requires only a scalar root-finding method.
  Such methods are relatively easy to implement from scratch, but also available in many open-source software packages.
  See also \Cref{rem:bisection}, below.
\medskip

  \noindent \textit{\textbf{Second,}} the transformation $\bm{\rho}_k \mapsto \sigma^{-1}(\bm{\rho}_k)$ in~\eqref{eq:PMD} is numerically unstable when $\bm{\rho}_k$ converges to a binary density field.
  Yet, a binary design is the desired outcome as $k \to \infty$.
  Directly tracking the latent variable $\bm{\psi}_k$ removes this instability.
  On the other hand, the original density variable can be reconstructed stably using the expression $\bm{\rho}_k=\sigma(\bm{\psi}_k)$ whenever necessary. 
\medskip

  \noindent \textit{\textbf{Third,}} working directly in the latent space provides a simple means to incorporate higher-order discretizations of the density variable, while still enforcing the box constraint $0 \leq \rho \leq 1$ at the discrete level.
  This aspect requires a more technical derivation of the method highlighted in the next subsection.
  See also \cite[Section~3.4]{simplmath}.

\paragraph*{Further remarks}
We close this subsection with a short list of remarks informed by our experiences using and teaching about the SiMPL method.

\begin{remark}[Computing the gradient]
  Even though~\eqref{eq:gradient-update} is written in terms of the latent variable $\bm{\psi}_k \in \mathbb{R}^{N_\rho}$, the gradient $\mathbf{g}_k$ (defined in~\eqref{eq:gradient}) continues to require differentiating $F$ with respect to $\bm{\rho}$.
  This gradient can be computed as usually done in topology optimization, using the adjoint method; cf.\ \Cref{lem:grad}, below.
\end{remark}

\begin{remark}[Quickly reaching binary designs]\label{rmk:density-latent}
  A binary design $\bm{\rho}^* = \lim_{k\to\infty} \sigma(\bm{\psi}_k)$ can be only achieved as the components of $\bm{\psi}_k$ approach $\pm\infty$.
  This is because each iteration of the design density $\bm{\rho}_k=\sigma(\bm{\psi}_k)$ satisfies the bound constraint strictly, i.e., $\mathbf{0} < \bm{\rho}_k < \mathbf{1}$.
  However, since $\sigma(-10) = \mathcal{O}(10^{-5})$ and $\sigma(-20)= \mathcal{O}(10^{-9})$, the SiMPL method achieves a sufficiently binary design once the components of $\bm{\psi}_k$ are on the order of $\pm 10$.
\end{remark}

\begin{remark}[Solving the volume projection equation]
\label{rem:bisection}
  To avoid possible numerical instabilities, we advocate for using a scalar root-finding method (e.g., the bisection method and the Illinois algorithm) to solve~\eqref{eq:vol-correction}.
  In this case, it is valuable to begin with well-defined upper and lower bounds on the solution.
  To this end, we assume that the previous iterate $\bm{\psi}_{k}$ satisfies the volume constraint $\mathbf{1}^\top\mathbf{M} \sigma(\bm{\psi}_{k})\leq\theta|\Omega|$.
  This inequality holds true for all $k \geq 1$ if $\bm{\psi}_{k}$ comes from a previous iteration of~\eqref{eq:full-update} and for $k = 0$ if we have started the SiMPL method with a feasible initial guess, such as $\bm{\psi}_{0} = \sigma^{-1}(\theta) \mathbf{1}$.
  In this case, owing to the monotonicity of $\sigma(x)=1/(1+\exp(x))$, we find that
  $$\mu_{k+1}\in [0,\max\{-\mathbf{g}_k\}],$$
  where the maximum is taken over all components of the negative gradient vector $-\mathbf{g}_k \in \mathbb{R}^{N_\rho}$.
  Thus, the root-finding method need only look within the interval $[0,\max\{-\mathbf{g}_k\}]$ for the root of~\eqref{eq:vol-correction}.
  In practice, we have found that the Illinois algorithm (a modified \textit{regula falsi} method, \cite{Dowell1971}) shows robust and fast convergence.
\end{remark}

\begin{remark}[Selecting the step sizes $\alpha_k$]
  The key to achieving efficiency with the SiMPL method is to use a (typically) \emph{increasing} sequence of step sizes $\alpha_k > 0$.
  Experience shows that setting $\alpha_k = \alpha_0 (k+1)^2$, where $\alpha_0>0$ is a tunable initial step size parameter, often leads to efficient solutions.
  However, we strongly advocate for using the line search strategies described in \Cref{sub:LineSearch} to achieve even better efficiency without parameter tuning.
\end{remark}

\begin{remark}[Convergence analysis]
  A rigorous convergence analysis of the SiMPL method at the function-space level can be found in the companion paper \cite{simplmath}.
\end{remark}

\begin{remark}[Taming the overflow]
  It is clear from~\Cref{rmk:density-latent} that $\sigma(x)$ converges to $0$ or $1$ exponentially as $x\rightarrow \pm\infty$. Therefore, we will obtain a numerically binary design when $\min\{|\bm{\psi}_k|\}$ is reasonably large.
  However, if we take a large number of steps or if step size becomes excessively large, then we may want to bound the latent variable to avoid numerical overflow.
  A straightforward approach is just projecting each component of the latent variable to the interval $[-M,M]$, where  $M\gg 0$ is some prescribed constant, i.e., $(\bm{\psi}_k)_i = \max\{\min\{(\bm{\psi}_k)_i, M\}, -M\}$, $i = 1,2, \ldots, N_\rho$.
  Another approach is regularizing the problem by adding an entropy penalty to the objective function~\eqref{eq:ObjectiveFunction}:
  \begin{equation}\label{eq:entropyPenalty}
      F(\bm{\rho})+\epsilon\varphi(\bm{\rho}),
  \end{equation}
  where $0 < \epsilon \ll 1$ is a small number.
  In this case, the latent variable update rule becomes
  \begin{equation*}
    \bm{\psi}_{k+1}=(1-\alpha_k\epsilon)\bm{\psi}_k-\alpha_k\mathbf{g}_k-\alpha_k\mu_{k+1}\mathbf{1}.
  \end{equation*}

\end{remark}

\subsection{Step size strategies}
\label{sub:LineSearch}

The key to achieving optimal efficiency with the SiMPL method is to use an \emph{increasing} sequence of step sizes $\alpha_k > 0$.
In this subsection, we outline some strategies for constructing such a sequence, strongly advocating for the line search strategy found in \Cref{alg:pmd-gbb}, below.
\begin{algorithm}[t]
  \caption{The SiMPL method}\label{alg:pmd-gbb}
  \algrenewcommand\algorithmicindent{0.75em}%
  \begin{algorithmic}[1]
    \Require exit tolerance $\mathtt{tol} > 0$ and $c_1 > 0$ (if using~\eqref{eq:armijo})
    \State $k \gets -1$
    \State$\bm{\rho}_{0} \gets \theta \mathbf{1}$
    \State $\bm{\psi}_0 \gets \sigma^{-1}(\bm{\rho}_{0})$
    \While{$\texttt{KKT}_k > \mathtt{tol}$} \Comment{Eq.~\eqref{eq:KKT_estimator}}
      \State $k \gets k+1$
      \State Evaluate the gradient $\mathbf{g}_k$ \Comment{Eq.~\eqref{eq:gradient}}
      \State $\displaystyle\alpha_k \gets \alpha_{k,0}$ \Comment{Eq.~\eqref{eq:stepsize}}
      \While{\texttt{true}}
      \State$\bm{\psi}_{k+1} \gets \bm{\psi}_k - \alpha_k (\mathbf{g}_k+\mu_{k+1})$ \Comment{Eq.~\eqref{eq:full-update}}
      \State$\bm{\rho}_{k+1} \gets \sigma(\bm{\psi}_{k+1})$ \Comment{Fig.~\ref{fig:mapping}}
      \If{\eqref{eq:linesearch} is satisfied}
      \State \textbf{break}
      \EndIf
      \State $\alpha_k \gets \alpha_k/2$
      \EndWhile
    \EndWhile

    \State \Return $\bm{\rho}_{k+1},\;F(\bm{\rho}_{k+1})$
  \end{algorithmic}
\end{algorithm}

\paragraph{Heuristics}
The SiMPL method was first tested in \cite{keith2023proximal} with a linearly growing step size, $\alpha_k = \alpha_0 (k+1)$, where $\alpha_0>0$ is a tunable initial step size parameter.
Further experience shows that the quadratic rule $\alpha_k = \alpha_0 (k+1)^2$ usually provides better efficiency.
In both cases, the number of iterations and the optimized design depend on the choice of the initial step size $\alpha_0 > 0$, demonstrating its influence on the efficiency and stability of the method.
Although there is value in verifying preliminary implementations with heuristics such as these, a robust and practical implementation should not rely on parameter tuning.
As a remedy, we advocate for the backtracking line search algorithm proposed in \cite{simplmath}.
This algorithm is able to adapt on-the-fly to problems with different scales and ensure the sequence of objective function values never increases, i.e., $F(\bm{\rho}_{k+1}) \leq F(\bm{\rho}_k)$ for every $k = 0, 1, 2, \ldots$

\paragraph{Backtracking line search}
Backtracking line search algorithms use function values to test a sufficient decrease condition and reduce the proposed step size if it fails.
For the SiMPL method, we propose two such conditions analyzed in \cite[Section~5]{simplmath}: the Armijo rule,
\begin{subequations}
  \label{eq:linesearch}
  \begin{equation}
    \label{eq:armijo}
      F(\bm{\rho}_{k+1})\leq F(\bm{\rho}_{k}) + c_1\mathbf{g}_k^\top\mathbf{M}(\bm{\rho}_{k+1}-\bm{\rho}_{k})
      ,
  \end{equation}
  where $0 < c_1 < 1$ is user-defined parameter, and the Bregman rule,
  \begin{multline}
    \label{eq:bregman}
      F(\bm{\rho}_{k+1}) \leq F(\bm{\rho}_{k}) +\mathbf{g}_k^\top\mathbf{M}(\bm{\rho}_{k+1}-\bm{\rho}_{k})\\+\frac{1}{\alpha_k}D_\varphi(\bm{\rho}_{k+1},\bm{\rho}_{k})
      .
  \end{multline}
\end{subequations}
Under relatively mild assumptions, both conditions ensure a monotonically-decreasing sequence of objective function values and guarantee convergence to a stationary point.
When $0 < c_1 \ll 1$, both conditions perform similarly, but the Bregman rule~\eqref{eq:bregman} has the apparent advantage that it is completely parameter-free.
On the other hand, we suggest using \eqref{eq:armijo} if the end user wants greater control over the line search selection process.
In this case, we recommend setting the default value of $c_1 = 10^{-4}$, and choosing larger values for a more conservative algorithm with shorter step sizes.

\paragraph{The Barzilai--Borwein step size}
Line search algorithms typically begin with a guess $\alpha_{k,0}$ for the next admissible step size $\alpha_{k}$.
We find that making this guess based on local curvature information significantly reduces the overall computational cost and accelerates convergence of the method.
Our starting point is the so-called long Barzilai--Borwein (BB) step size \cite{barzilai1988bbm}, which can be viewed as an approximated local Lipschitz continuity constant.
It is defined via the previous two iterates $\bm{\rho}_{k},\bm{\rho}_{k-1}$ and search directions $\mathbf{g}_k,\mathbf{g}_{k-1}$:
\begin{equation*}
  \alpha_{k, \rm BB}
  =
  \frac{(\bm{\rho}_{k}-\bm{\rho}_{k-1})^\top\mathbf{M}(\bm{\rho}_{k}-\bm{\rho}_{k-1})}{|(\mathbf{g}_k-\mathbf{g}_{k-1})^\top\mathbf{M}(\bm{\rho}_{k}-\bm{\rho}_{k-1})|}
  \,.
\end{equation*}
Here, the absolute value is used in the denominator to ensure the positivity of $\alpha_{k, \rm BB}$.

\paragraph{Generalizing the BB step size}
For the mirror descent method, we choose to generalize the BB step size as follows:
\begin{equation}
\label{eq:GBB}
    \alpha_{k, \rm GBB}  = \frac{(\bm{\psi}_k-\bm{\psi}_{k-1})^\top\mathbf{M}(\bm{\rho}_{k}-\bm{\rho}_{k-1})}{|(\mathbf{g}_k-\mathbf{g}_{k-1})^\top\mathbf{M}(\bm{\rho}_{k}-\bm{\rho}_{k-1})|}
    \,.
\end{equation}
In this case, the step size $\alpha_{k, \rm GBB}$ can be understood as an approximation of the relative continuity constant, a generalization of the Lipschitz constant that often plays an important role in the convergence analysis of mirror descent methods \cite{lu2019relative,teboulle2018nolip,bauschke2017-NoLip}; see also \cite[Section~5.1]{simplmath}.

\paragraph{Estimating the next step size}
We often find that~\eqref{eq:GBB} selects an exponentially growing step size guess.
To help avoid over-estimates, we suggest taking the geometric mean of $\alpha_{k, \rm GBB}$ with the previous step size as the line search step size guess:
\begin{subequations}
\label{eq:stepsize}
\begin{equation}
\label{eq:stepsize1}
  \alpha_{k,0}=\sqrt{\alpha_{k, \rm GBB}\alpha_{k-1}}
  \,,
\end{equation}
for each $k = 1, 2, \ldots$
Clearly, previous information is not available at first iteration, $k = 0$, thus we set
\begin{equation}
\label{eq:stepsize2}
  \alpha_{0,0} = 1/\max\{|\mathbf{g}_k|\}
  \,.
\end{equation}
\end{subequations}

\subsection{Stopping criteria} 
\label{sub:stopping_criteria}

The SiMPL method, with the line search strategy introduced in \Cref{sub:LineSearch}, often decreases objective function values rapidly.
However, important design changes can occur when the function value increments, $\delta F_k = F(\bm{\rho}_{k}) - F(\bm{\rho}_{k+1})$, are very small. 
Moreover, since the size of these increments is influenced by the step sizes $\alpha_k$, we do not recommend relying on sufficiently small $\delta F_k$ as the only stopping criterion.
At the very least, we advocate for also estimating the KKT stationary condition.

\paragraph{KKT conditions}
The Karush--Kuhn--Tucker (KKT) conditions \cite{nocedal1999numerical,beck2023introduction} are necessary (and sometimes sufficient) optimality conditions for a solution of a constrained optimization problem.
Denote a local minimizer of~\eqref{eq:ReducedProblem} by $\bm{\rho}^\star \in \mathcal{A}_h$ and the gradient of $F$ at $\bm{\rho}^{\star}$ by $\mathbf{g}^{\star} = \mathbf{M}^{-1}\mathbf{d} F(\bm{\rho}^{\star})$. 
The KKT conditions imply the existence of Lagrange multipliers $\mu^{\star} \geq 0$ and $\bm{\lambda}^{\star} \in \mathbb{R}^{N_\rho}$ satisfying the stationarity equation
\begin{subequations}
\label{eq:KKTConditions}
\begin{equation}
\label{eq:stationarity}
  \mathbf{g}^{\star} + \bm{\lambda}^{\star} + \mu^{\star}\mathbf{1}
  =
  \mathbf{0}
  \,,
\end{equation}
together with the complementarity conditions
\begin{equation}
\label{eq:complementarity_mu}
  \mu^{\star} = 0
  ~~\text{ if }
  \mathbf{1}^\top\mathbf{M}\bm{\rho}<\theta|\Omega| 
\end{equation}
and
\begin{equation}
\label{eq:lambda-signs}
  (\bm{\lambda}^{\star})_i
  \begin{cases}
    \geq {0} & \text{if } (\bm{\rho}^{\star})_i = {1}, \\
    \leq {0} & \text{if } (\bm{\rho}^{\star})_i = {0}, \\
    = {0} & \text{if } {0} < (\bm{\rho}^{\star})_i < {1},
  \end{cases}
\end{equation}
for each $1\leq i \leq N_\rho$.
\end{subequations}

\paragraph{Approximate Lagrange multiplier}
We suggest stopping the SiMPL method when the iterate $\bm{\rho}_{k+1}$ satisfies~\eqref{eq:KKTConditions} sufficiently accurately.
To this end, we manipulate~\eqref{eq:full-update} to derive the following identity:
\[
  \mathbf{g}_k + \frac{{\bm{\psi}}_{k+1}-\bm{\psi}_k}{\alpha_k} + {\mu}_{k+1}\mathbf{1}
  =
  \mathbf{0}
  \,,
\]
where ${\mu}_{k+1} \geq 0$ satisfies the complementarity condition~\eqref{eq:complementarity_mu} by construction.
The similarity to~\eqref{eq:stationarity} is not coincidental \cite[Proposition~4.8]{simplmath}, and suggests defining
\begin{equation}
\label{eq:ApproximateLM}
  \bm{\lambda}_k := ({\bm{\psi}}_{k+1}-\bm{\psi}_k)/{\alpha_k}
\end{equation}
as an approximation to the Lagrange multiplier $\bm{\lambda}^{\star}$.

\paragraph{KKT estimator}
It remains to measure how well $\bm{\lambda}_k$ satisfies the inequalities in~\eqref{eq:lambda-signs}.
With this goal in mind, we introduce the following positive vector $\bm{\eta}_k \in \mathbb{R}^{N_\rho}_+$ encoding the component-wise violations of complementarity condition~\eqref{eq:lambda-signs}:
\begin{subequations}
\label{eq:KKT_estimator}
\begin{equation}
  \label{eq:kkt-brendan}
  \bm{\eta}_k  =\max\{-\bm{\rho}_k\bm{\lambda}_k, (\mathbf{1}-\bm{\rho}_k)\bm{\lambda}_k\}.
\end{equation}
However, we note that there is no unique definition of $\bm{\eta}_k$, and we have also found promising results (cf.\ \cite[Section~6]{simplmath}) with the alternative choice
\begin{equation}
  \label{eq:kkt-thomas}
    \bm{\eta}_k  = \bm{\lambda}_k - \min\{\mathbf{0}, \bm{\rho}_k + \bm{\lambda}_k\}\\ - \max\{\mathbf{0}, \bm{\rho}_k - \mathbf{1} + \bm{\lambda}_k\}
    \,.
\end{equation}
In either case, once $\bm{\eta}_k$ has been specified, we suggest stopping the SiMPL method once
\begin{equation}
  \texttt{KKT}_k := \mathbf{1}^\top \mathbf{M} \bm{\eta}_k
  \leq \mathtt{tol}
  \,,
\end{equation}
where $\mathtt{tol}> 0$ is a prescribed accuracy tolerance; cf.\ line 4 of \Cref{alg:pmd-gbb}.
\end{subequations}

\subsection{First optimize then discretize}\label{sub:opt-disc}
\begin{subequations}
\label{eq:conti-prob}
In this subsection, we rederive the SiMPL method using the optimize-then-discretize paradigm.
In this case, once the method is established at the function space level, we show how it may be used to derive high-order discrete SiMPL methods for more general types of meshes.
When the lowest-order discretization is used for $\rho_h\in Q_h$, the resulting method is equivalent to the one derived in \Cref{sub:FDTO}.
However, using  the optimize-then-discretize paradigm, we can also derive the SiMPL method for high-order discretizations of the density variable without losing feasibility.
Finally, the resulting formulation shows mesh- and degree-independent behavior, see \Cref{fig:mbb-mesh} and \cite{simplmath}.

\paragraph*{Problem definition}
Consider the following topology optimization problem analogous to \eqref{eq:disc-prob} but formulated in function spaces:
  \begin{gather}
  \label{eq:conti-objective}
    \operatorname{minimize}\ \widehat{F}(\tilde{\rho},u)
  \end{gather}
  over $\rho \in L^2(\Omega)$, $\tilde{\rho} \in H^1(\Omega)$, and $u \in V\subset [H^1(\Omega)]^d$ with $d = 2$ or $3$, subject to
  \begin{gather}
    \label{eq:conti-state-eq}          \int_\Omega \big(r(\tilde{\rho})\mathsf{C}\,\varepsilon(u)\big) : \varepsilon(v) \dd x = \int_\Omega f\cdot v\dd x,
    \\
    \label{eq:conti-filter-eq}         \int_\Omega \epsilon^2\nabla \tilde{\rho}\cdot\nabla \tilde{q}+\tilde{\rho}\tilde{q}\dd x=\int_\Omega \rho \tilde{q} \dd x,
    \intertext{for all $v\in V$ and $\tilde{q}\in H^1(\Omega)$ and}
    \label{eq:conti-bound-constraint}  0\leq\rho(x)\leq1~\text{ for almost every } x \in \Omega,
    \\
    \label{eq:conti-volume-constraint} \int_\Omega \rho\dd x\leq\theta|\Omega|.
  \end{gather}
\end{subequations}
Here, $\mathsf{C}$ denotes the (fourth-order) elasticity tensor, $\varepsilon(u) = (\nabla u +\nabla u^\top)/2$ denotes the symmetric gradient, $r(\tilde{\rho}) = \rho_0 + \tilde{\rho}^p (1 - \rho_0)$, with exponent $p> 1$ and nominal density $0 < \rho_0 \ll 1$, denotes the continuous form of the SIMP penalization law, and $f \in L^2(\Omega)$ denotes an applied load.
All other parameters are the same as in \Cref{sub:FDTO}.
Following \cite{simplmath}, we define the set of admissible density functions to be
\begin{multline*}
  \mathcal{A}
  =
  \bigg\{
    \rho\in L^2(\Omega) \mid \int_\Omega \rho \dd x \leq \theta |\Omega| \text{ and }\\
    0\leq \rho(x)\leq 1 \text{ for almost every } x \in \Omega
  \bigg\}
  \,.
\end{multline*}
Then, assuming the objective function in~\eqref{eq:conti-objective} is sufficiently smooth, we rewrite problem~\eqref{eq:conti-prob} using a so-called reduced objective function, written solely as a function of the density $\rho$.
In particular, we write
\begin{equation}\label{eq:red-top-opt}
  \min_{\rho \in \mathcal{A}} F(\rho)
  \,,
\end{equation}
where $F(\rho) := \widehat{F}(\tilde{\rho}(\rho),u(\tilde{\rho}(\rho)))$.

\paragraph{The gradient}
Before we derive the SiMPL method in infinite-dimensional function spaces, we introduce a result from \cite{simplmath} with a derivation given in \Cref{sec:appendix}.
In what follows, $F^\prime(\rho)$ denotes the Fr\'echet derivative of $F$ at $\rho$ and $\langle\cdot,\cdot\rangle$ is the natural duality pairing on the implied function spaces.
\vspace*{-1em}

\begin{proposition}\label{lem:grad}
  Given a design density $\rho \in \mathcal{A}$, let ${u}\in V$ and $\tilde{\rho}\in H^1(\Omega)$ be the unique solutions to \eqref{eq:conti-state-eq} and \eqref{eq:conti-filter-eq}, respectively.
  If we assume that $\widehat{F}$ is continuously Fr\'{e}chet differentiable, then the reduced objective function $F$ in~\eqref{eq:red-top-opt} is Fr\'{e}chet differentiable in $L^\infty(\Omega)$ and its Fr\'{e}chet derivative at $\rho$ can be obtained by solving the following sequence of adjoint problems:
  Find $\lambda\in V$ such that
  \begin{subequations}
  \label{eqs:grad}
  \begin{align}
  \label{eq:lambda}
    \int_\Omega \big(r(\tilde{\rho})\mathsf{C}\,\varepsilon(\lambda)\big):\varepsilon(v)\dd x= \langle\partial_u\widehat{F}(\tilde{\rho},u), v \rangle
  \end{align}
  for all $v \in V$ and then find $\tilde{g} \in H^1(\Omega)$ such that
  \begin{multline}
  \label{eq:gtilde}
     \int_\Omega \epsilon^2\nabla\tilde{g}\cdot\nabla \tilde{q}+\tilde{g}\tilde{q}\dd x
     =\langle \partial_{\tilde{\rho}} \widehat{F}(\tilde{\rho},u), \tilde{q} \rangle\\
      -\int_\Omega \big(r'(\tilde{\rho})\mathsf{C}\,\varepsilon(u):\varepsilon(\lambda)\big)\tilde{q}\dd x
  \end{multline}
  \end{subequations}
  for all $\tilde{q} \in H^1(\Omega)$.
  In particular, we have that
  \begin{equation}\label{eq:derivative}
    \langle F'(\rho), q \rangle =
    \int_\Omega \tilde{g}q\dd x ~\text{ for all } q\in L^\infty(\Omega).
  \end{equation}
\end{proposition}
\vspace*{-0.5em}

\noindent
We refer to $\tilde{g}$ in~\eqref{eq:derivative} as the \emph{gradient} of $F$ at $\rho$.

\paragraph{The continuous SiMPL method}
We are now ready to define the local energy functional $J_\varphi(\rho;\rho_k)$ that is minimized at each iteration of the SiMPL method:
\begin{equation}
\label{eq:CtsSubproblem}
  \rho_{k+1} = \argmin_{\rho \in \mathcal{A}} J_\varphi(\rho;\rho_k)
  \,.
\end{equation}
In particular, we define
\begin{equation*}
  J_\varphi(\rho;\rho_k)
  =
  \int_\Omega \tilde{g}_k \rho 
  +\frac{1}{\alpha_k}D_\varphi(\rho,\rho_k) \dd x,
\end{equation*}
where $\tilde{g}_k$ is the gradient of $F$ at $\rho_k$ and
\begin{subequations}\label{eq:conti-varphi}
\begin{equation*}
  D_\varphi(\rho,q)
  =
  \varphi(\rho)-\varphi(q) + \varphi^\prime(q)(\rho-q)
  \,,
\end{equation*}
is the Bregman divergence associated to the Fermi--Dirac entropy,
\begin{equation*}
  \varphi(\rho) = \int_\Omega \rho\ln(\rho)+(1-\rho)\ln(1-\rho)\dd x
  \,.
\end{equation*}
\end{subequations}
\Cref{eq:CtsSubproblem} is analyzed rigorously in \cite[Theorem~3.4]{simplmath}, revealing the following update formula:
\[
 \rho_{k+1} = \sigma\big(\sigma^{-1}(\rho_k)-\alpha_k \tilde{g}_k-\alpha_k\mu_{k+1}\big)
 \,.
\]
We then derive the following two-stage formulae by introducing the latent variable $\psi_k = \sigma^{-1}(\rho_k)$:
\begin{subequations}\label{eq:full-update-infinite}
  \begin{align}
    \label{eq:gradient-update-infinite}
    \psi_{k+1/2}        & =\psi_{k}-\alpha_k \tilde{g}_k,                          \\
    \label{eq:volume-update-infinite}
    \psi_{k+1}          & =\psi_{k+1/2}-\alpha_k\mu_{k+1}.
  \end{align}
  where $\mu_{k+1} \geq 0$ solves the non-smooth volume correction equation
  \begin{multline}
    \label{eq:vol-correction-infinite}
     \min\bigg\{
      \mu_{k+1},\\
      \theta|\Omega| - \int_\Omega \sigma( \psi_{k+1/2}-\alpha_k\mu_{k+1})\dd x
     \bigg\}
     = 0
    \,.
  \end{multline}
  In particular, $\mu_{k+1}= 0$ when $\int_\Omega \sigma(\psi_{k+1/2}-\alpha_k\mu_{k+1})\dd x < \theta|\Omega|$. Otherwise, $\mu_{k+1} \geq 0$ is the unique non-negative number satisfying
\end{subequations}
  \begin{equation*}
    \label{eq:vol-correction-infinite2}
    \int_\Omega \sigma( \psi_{k+1/2}-\alpha_k\mu_{k+1})\dd x = \theta|\Omega|.
  \end{equation*}
Even though $\psi_k$ is not expected to be bounded as $k \to \infty$, analysis shows that $L^\infty(\Omega)$ is the natural function space for the latent variables in~\eqref{eq:full-update-infinite}; cf.\ \cite[Section~6.3]{keith2023proximal}.
Fortunately, under mild assumptions on the domain $\Omega$ and the external load $f$ \cite[Section~2.2]{simplmath}, we can guarantee that each $\tilde{g}_k$ belongs to $L^\infty(\Omega)$.
Thus,~\eqref{eq:full-update-infinite} is always well-defined so long as the algorithm begins at a feasible initial guess $\psi_0 \in L^\infty(\Omega)$.

\paragraph{Discretizing the SiMPL method}
Let $Q_h \subset L^\infty(\Omega)$ be a finite element subspace with ordered basis $\bm{\Phi} = (\phi_1,\phi_2,\ldots,\phi_{N_\rho})^\top$, where each $\phi_i \in L^\infty(\Omega)$.
Discretizing~\eqref{eq:full-update-infinite} with finite elements requires expressing the approximations to each $\psi_k$ and $\tilde{g}_k$ as linear combinations of these basis functions.
In particular, we write
\begin{equation}
\label{eq:SiMPL_discrete_variables}
  \psi_k \approx \bm{\psi}_k^\top \bm{\Phi}
  ~~\text{ and }~~
  \tilde{g}_k \approx \mathbf{g}_k^\top \bm{\Phi}
  \,,
\end{equation}
where $\bm{\psi}_k, \mathbf{g}_k \in \mathbb{R}^{N_{\rho}}$ are coefficient vectors.
If $\bm{\Phi}$ forms a partition of unity, i.e., $\sum \phi_i(x) = 1$ for all $x\in\Omega$, then we uncover the following discretized algorithm:
\begin{align*}
  \bm{\psi}_{k+1/2} & =\bm{\psi}_{k}-\alpha_k\mathbf{g}_k,
  \\
  \bm{\psi}_{k+1}   & =\bm{\psi}_{k+1/2}-\alpha_k\mu_{k+1}\mathbf{1}
  ,
\end{align*}
where, similar to before, $\mu_{k+1} \geq 0$ comes from solving
\begin{equation*}
  \int_\Omega \sigma( \bm{\psi}_{k+1/2}^\top \bm{\Phi}-\alpha_k\mu_{k+1})\dd x = \theta|\Omega|
\end{equation*}
if the volume constraint $\int_\Omega \sigma( \bm{\psi}_{k+1/2}^\top \bm{\Phi}-\alpha_k\mu_{k+1})\dd x \leq \theta|\Omega|$ is violated.
This algorithm coincides with~\eqref{eq:full-update} when $\bm{\Phi}$ is composed of the piecewise-constant indicator functions for cells in a grid so long as the gradient $\tilde{g}_k$ is discretized following~\eqref{eq:derivative}.

\paragraph{Discretizing the gradient}
\Cref{eq:conti-state-eq,eq:conti-filter-eq,eq:lambda,eq:gtilde} must be discretized and solved at each iteration of the SiMPL method to generate an approximation of $\tilde{g}_k$.
In particular, solving~\eqref{eq:gtilde} will return an approximation
\begin{align}
\label{eq:SiMPL_discrete_gradient}
  \tilde{g}_k \approx \tilde{\mathbf{g}}_k^\top \tilde{\bm{\Phi}}
\end{align}
belonging to $\tilde{Q}_h \subset H^1(\Omega)$.
Here, $\tilde{\bm{\Phi}} = (\tilde{\phi}_1, \tilde{\phi}_2,\ldots, \tilde{\phi}_{N_{\tilde{\rho}}})^\top$ is an ordered basis satisfying $\tilde{\phi}_j \in H^1(\Omega)$.
Since the subspaces $Q_h$ and $\tilde{Q}_h$ used in~\eqref{eq:SiMPL_discrete_variables} and \eqref{eq:SiMPL_discrete_gradient}, respectively, will generally not coincide, we require a formula relating the coefficient vectors ${\mathbf{g}}_k$ and $\tilde{\mathbf{g}}_k$.
The most natural formula arises by defining $\mathbf{g}_k^\top \bm{\Phi}$ to be the Galerkin projection of $\tilde{\mathbf{g}}_k^\top \tilde{\bm{\Phi}}$ onto $Q_h$; i.e., by finding the unique $\mathbf{g}_k \in \mathbb{R}^{N_\rho}$ that satisfies
\begin{equation}
\label{eq:GalerkinProjection}
  \int_\Omega \mathbf{g}_k^\top \tilde{\bm{\Phi}} \phi_i \dd x
  =
  \int_\Omega \tilde{\mathbf{g}}_k^\top \tilde{\bm{\Phi}} \phi_i \dd x
  \,,
\end{equation}
for each $i = 1,2,\ldots, N_\rho$.
Equivalently, we may follow~\eqref{eq:derivative}, to rewrite~\eqref{eq:GalerkinProjection} as a linear equation relating $\mathbf{g}_k$ to the coefficient vector representation of the Fr\'echet derivative of $F$: namely,
\begin{equation}
\label{eq:DiscreteGradientEquation}
  \mathbf{M}\mathbf{g}_k = \mathbf{d}F_k
  \,,
\end{equation}
where $\mathbf{d}F_k := \mathbf{N}\tilde{\mathbf{g}}_k$ and $\mathbf{N}_{ij} = \int_\Omega \phi_i\tilde{\phi}_j \dd x$ is a (non-symmetric) mass matrix.

\paragraph{High-order discretizations}

An appealing feature of the SiMPL method is that it guarantees bound-preserving discrete design densities.
Indeed, no matter the form of the basis functions used to approximate $\psi_k$ in~\eqref{eq:SiMPL_discrete_variables}, the sigmoid function returns a bound-preserving density
\[
  \sigma(\bm{\psi}_k^\top \bm{\Phi}) \approx \rho_k
  \,.
\]
Clearly, $0 \leq \sigma(\bm{\psi}_k^\top \bm{\Phi}) \leq 1$ by construction.
Numerical experiments with the SiMPL method supporting the use of high-order discrete densities can be found in \cite{simplmath} and Problem 3 in the next section.

\begin{figure*}[t]
  \centering
  \begin{tabular}{c|c|c|c|c|}
    $k$                                                                       & \textbf{SiMPL-A} & \textbf{SiMPL-B} & OC & MMA \\
    10                                                                        &
    \includegraphics[width=0.21\textwidth]{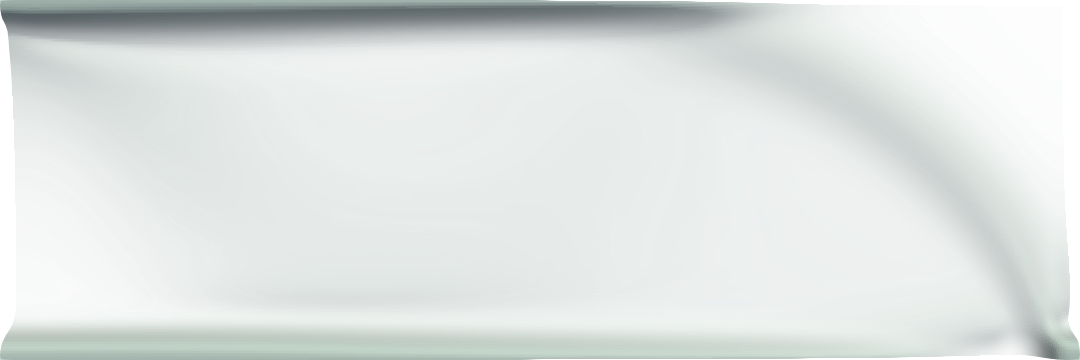} &
    \includegraphics[width=0.21\textwidth]{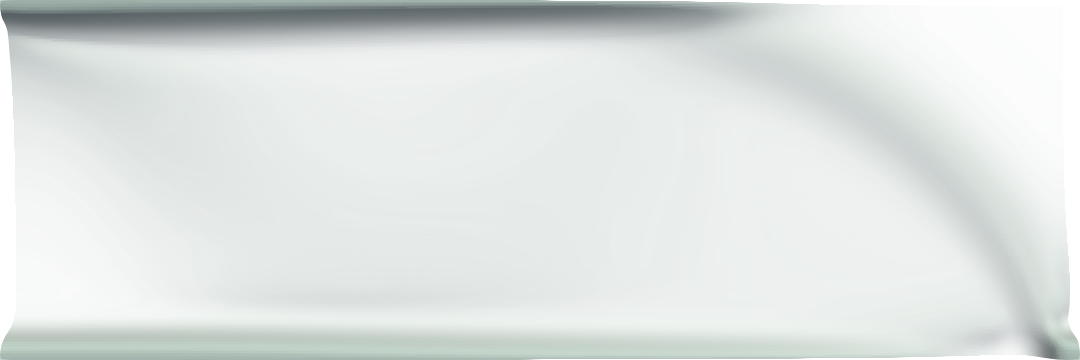} &
    \includegraphics[width=0.21\textwidth]{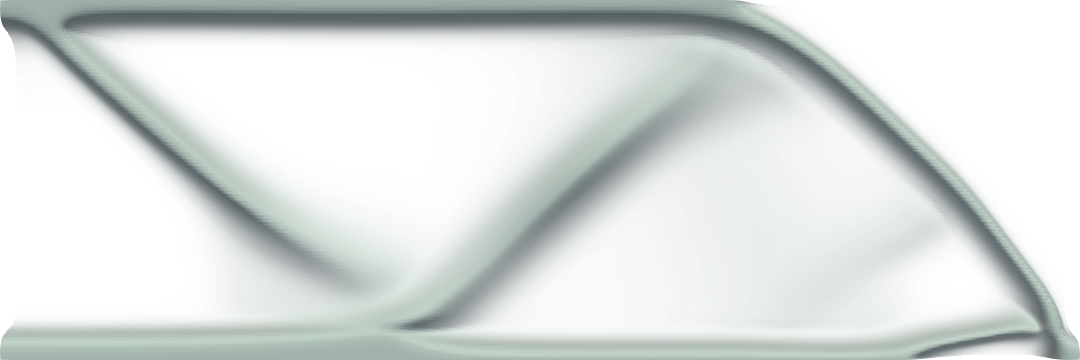}           &
    \includegraphics[width=0.21\textwidth]{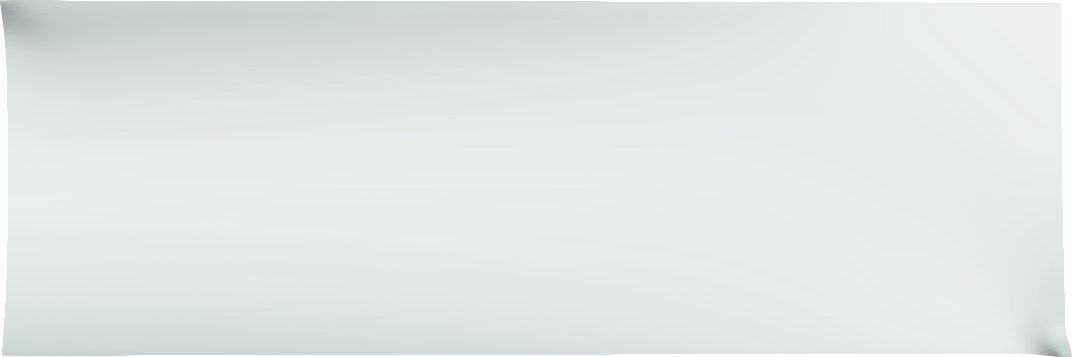}                                                          \\
    20                                                                        &
    \includegraphics[width=0.21\textwidth]{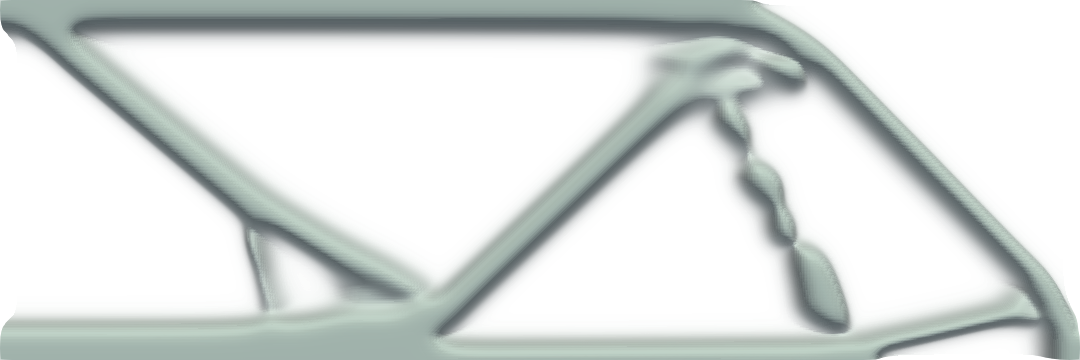} &
    \includegraphics[width=0.21\textwidth]{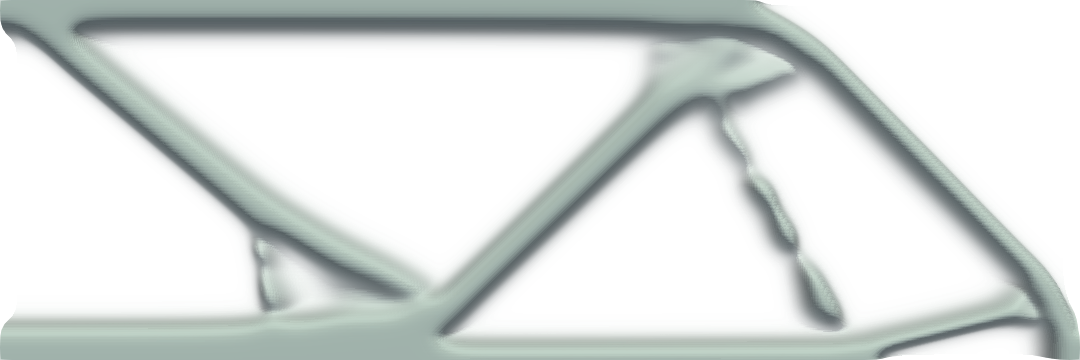} &
    \includegraphics[width=0.21\textwidth]{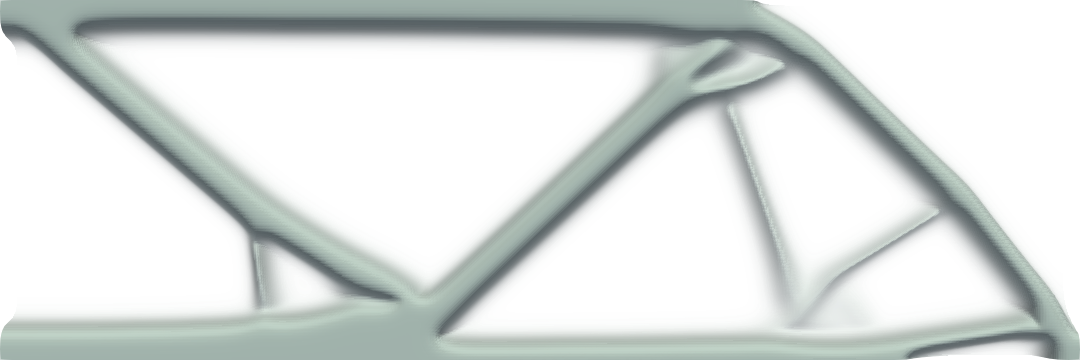}           &
    \includegraphics[width=0.21\textwidth]{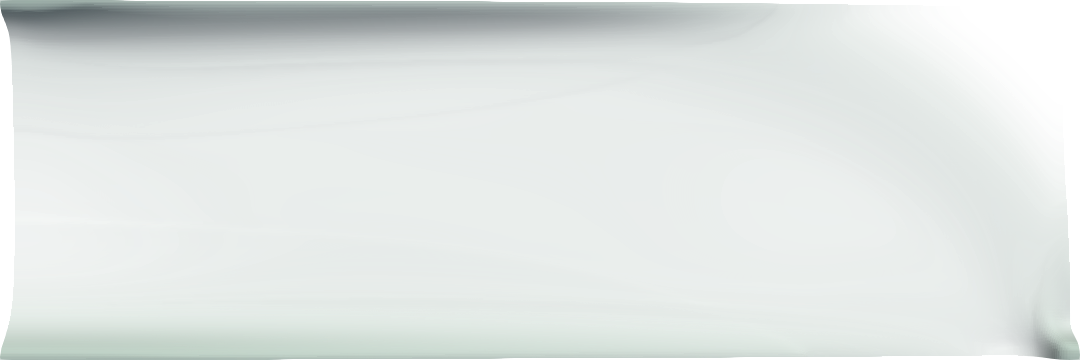}                                                          \\
    30                                                                        &
    \includegraphics[width=0.21\textwidth]{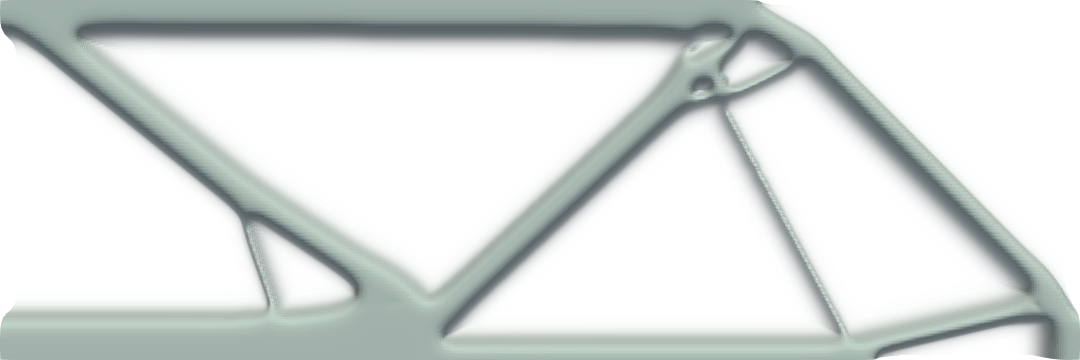} &
    \includegraphics[width=0.21\textwidth]{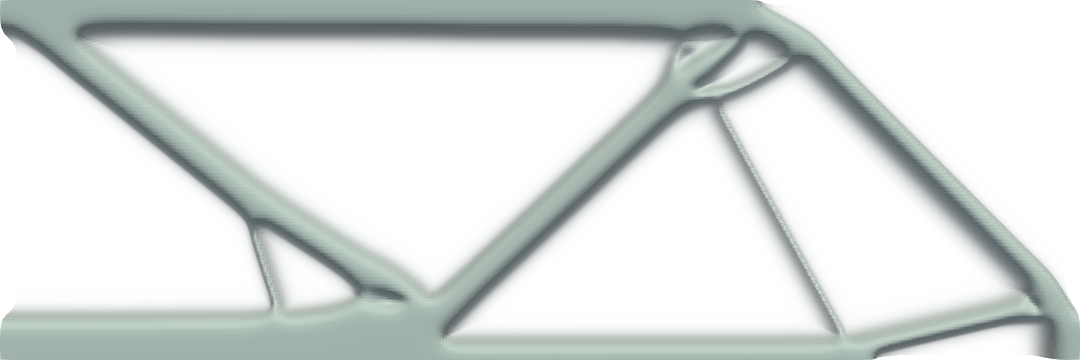} &
    \includegraphics[width=0.21\textwidth]{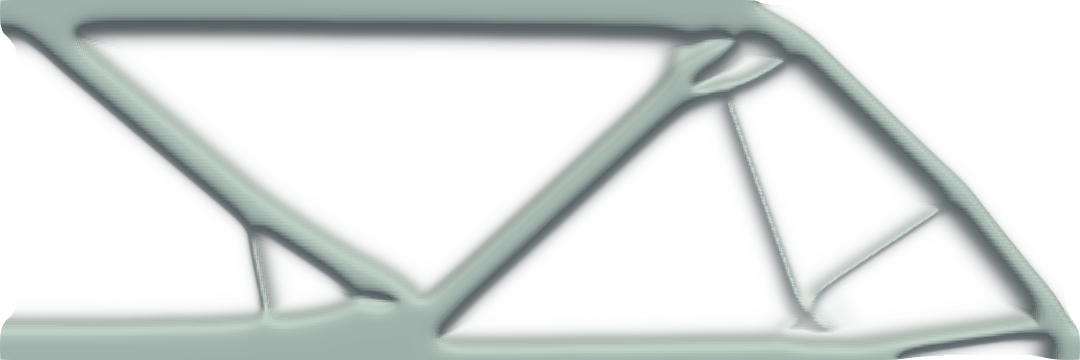}           &
    \includegraphics[width=0.21\textwidth]{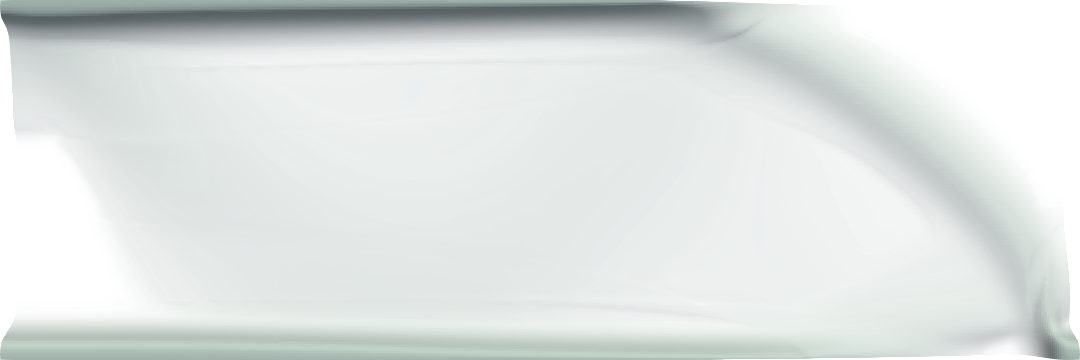}                                                          \\
    50                                                                        &
    \includegraphics[width=0.21\textwidth]{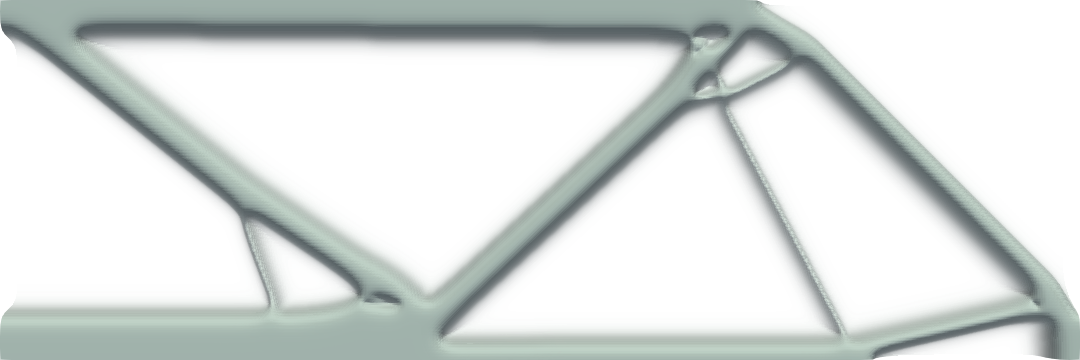} &
    \includegraphics[width=0.21\textwidth]{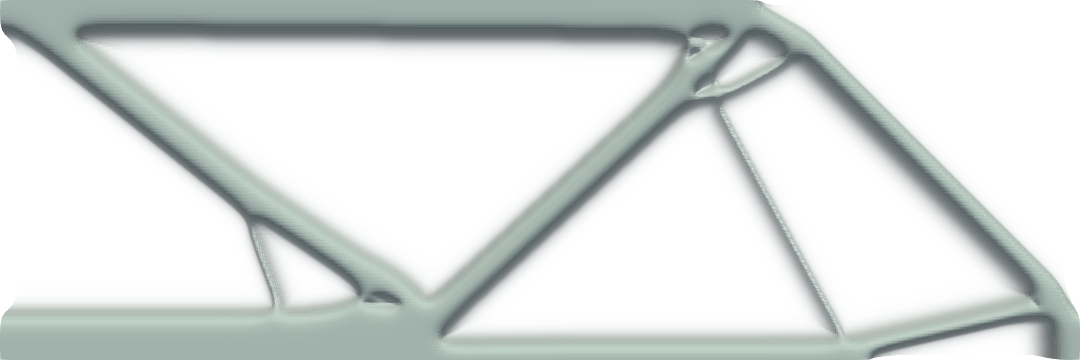} &
    \includegraphics[width=0.21\textwidth]{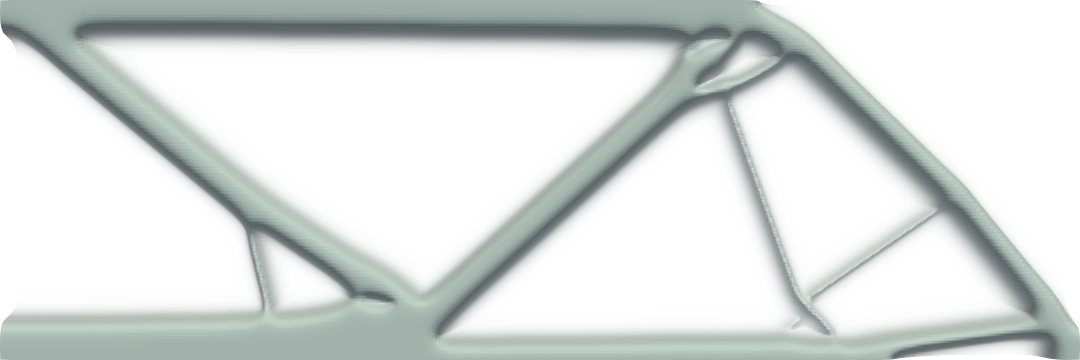}           &
    \includegraphics[width=0.21\textwidth]{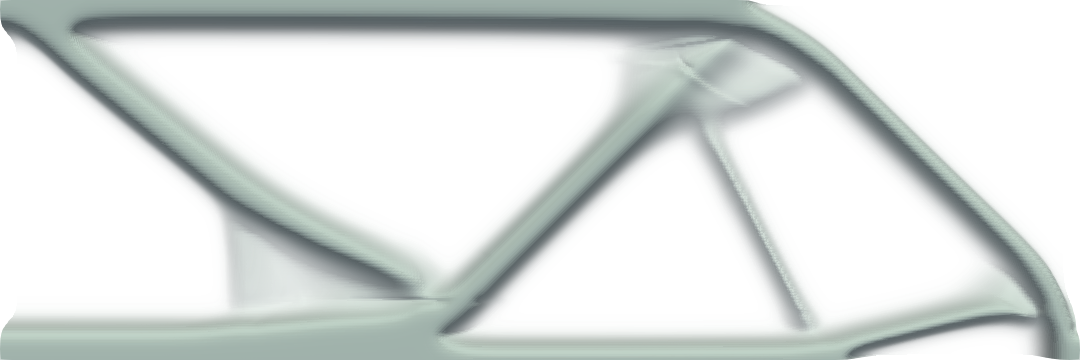}                                                          \\
    100                                                                       &
                                                                              &                  &
    \includegraphics[width=0.21\textwidth]{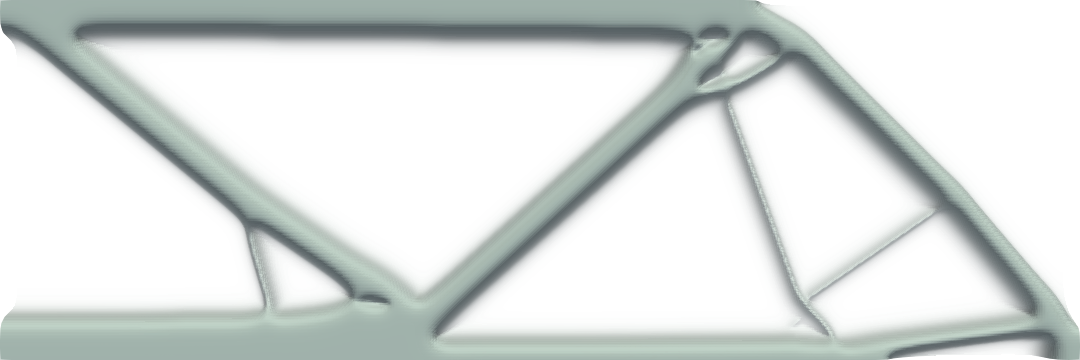}           &
    \includegraphics[width=0.21\textwidth]{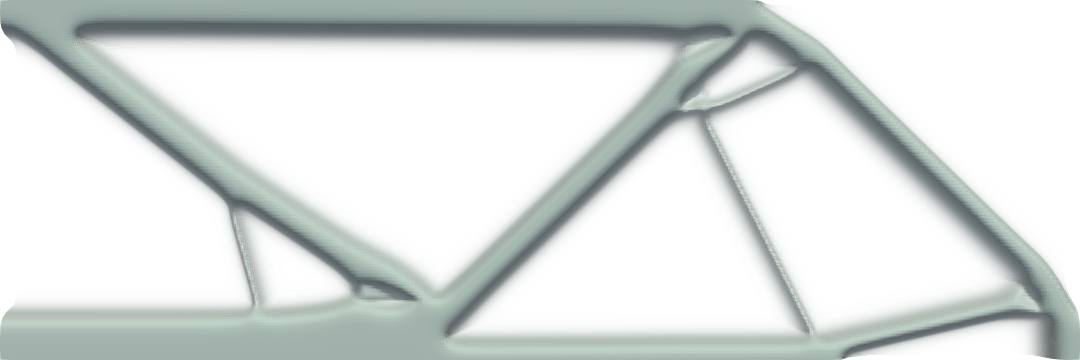}                                                          \\
    300                                                                       &
                                                                              &                  &
    \includegraphics[width=0.21\textwidth]{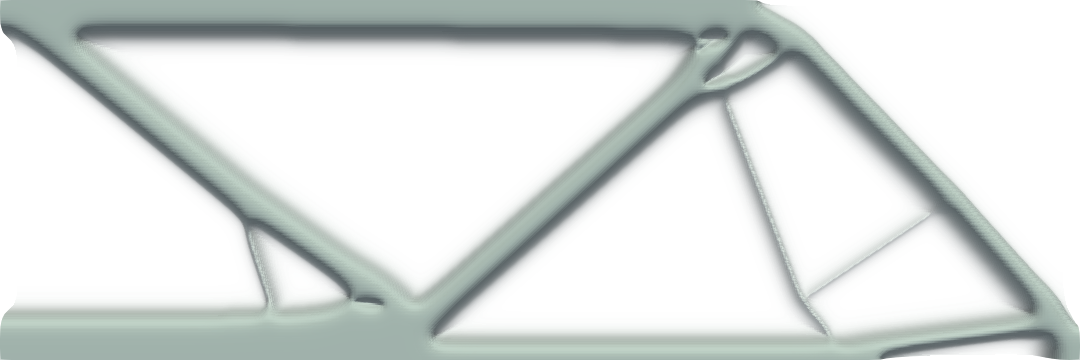}           &
    \includegraphics[width=0.21\textwidth]{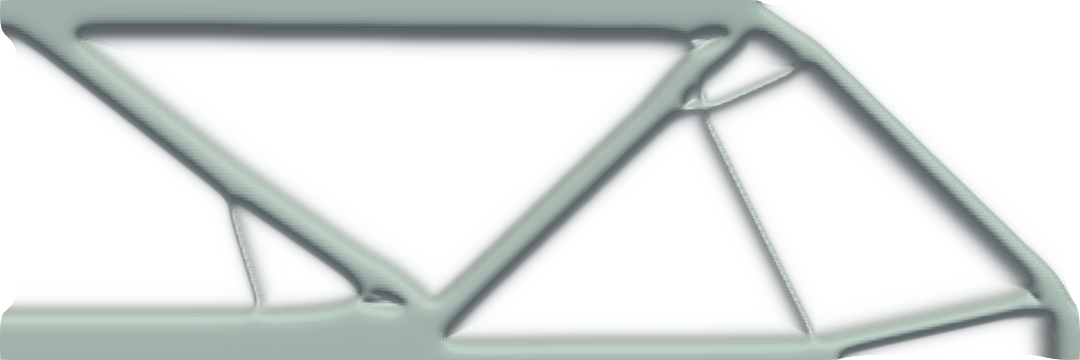}                                                          \\
  \end{tabular}
  \caption{Problem 1. Filtered density, $\tilde{\rho}$, for selected iterations $k$. From left to right: SiMPL-A, SiMPL-B, OC, and MMA. The final number of iterations are 50 (SiMPL-A), 46 (SiMPL-B), 300 (OC), and 300 (MMA).}
  \label{fig:mbb}
\end{figure*}

\section{Applications}\label{sec:app}
In this section, we report our findings from applying the SiMPL method to several TO problems.
We also compare SiMPL to the popular OC and MMA algorithms.

\paragraph*{Set-up}
For simplicity and consistency, we always follow the finite element discretization used in~\eqref{eq:LowestOrderSpaces} and \eqref{eq:disc-prob}.
To demonstrate the flexibility of the method, we use two versions of \Cref{alg:pmd-gbb}, which we refer to as SiMPL-A and SiMPL-B.
SiMPL-A uses the Armijo rule \eqref{eq:armijo}, while SiMPL-B uses the Bregman rule~\eqref{eq:bregman}.
We always set $c_1=10^{-4}$ in \eqref{eq:armijo} for SiMPL-A.
Unless otherwise specified, we always set $\bm{\rho}_0=\theta \mathbf{1}$ as the initial design density, use the complementarity vector $\bm{\eta}_k$ defined in~\eqref{eq:kkt-brendan}, and fix $r_{\rm min} = 0.02$ as the filter radius, implying $\epsilon=0.02/(2\sqrt{3})$ in~\eqref{eq:filter-eq}; cf.\ \cite{lazarov2011-filter}.
Finally, the implemented OC algorithm follows the formulation in \cite{Andreassen2011}, and MMA is executed with default parameters given in \cite{aage2015petsc}.
Both the OC and MMA updates are limited to obtain stable convergence behavior,
  \begin{equation*}
    \ell_i=\max\{0,(\bm{\rho}_k)_i-\texttt{ch}\},~ u_i=\min\{1,(\bm{\rho}_k)_i+\texttt{ch}\}.
  \end{equation*}
  Here, all results are reported with \texttt{ch}=0.15 which gives the best performance among the tested values \texttt{ch} $\in\{0.05,0.1,0.15,0.2,0.25,0.3,0.4\}$.
  The results of MMA can be further improved by tuning the internal parameters, but we do not pursue this here.

\subsection{Compliance minimization}
The first application is compliance minimization:
\begin{equation*}
  \begin{aligned}
    \min_{\bm{\rho}}\  & \mathbf{f}^\top \mathbf{u}
    \\
    \text{subject to }~
                       &
                       \mathbf{K}(\tilde{\bm{\rho}})\mathbf{u} =\mathbf{f}                                                   \\
                       & (\epsilon^2\mathbf{A}+\tilde{\mathbf{M}})\tilde{\bm{\rho}} =\mathbf{N}\bm{\rho}  \\
                       & \mathbf{0}\leq \bm{\rho} \leq \mathbf{1},                                        \\
                       & \mathbf{1}^\top\mathbf{M}\bm{\rho} \leq\theta|\Omega|
                       \,,
  \end{aligned}
\end{equation*}
where $\mathbf{f}$ is the external force.
In particular, we first consider the popular MBB (Messerschmitt--B{\"o}lkow--Blohm) beam problem \cite{Bendsoe2004,Andreassen2011}.

\paragraph*{Problem 1: 2D MBB beam}
\begin{figure*}
	\centering
    ~
  \includegraphics[width=0.36\textwidth]{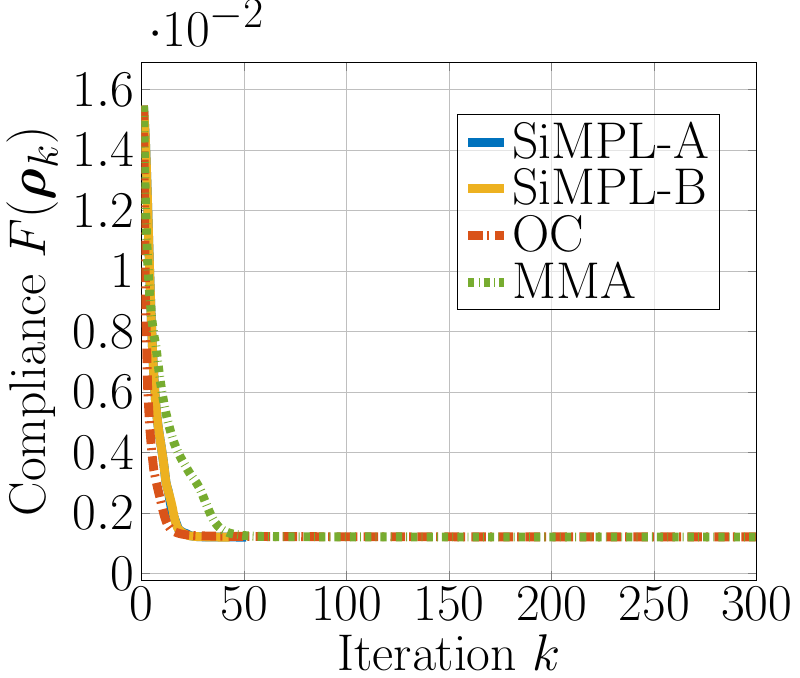}
  ~
  \includegraphics[width=0.38\textwidth]{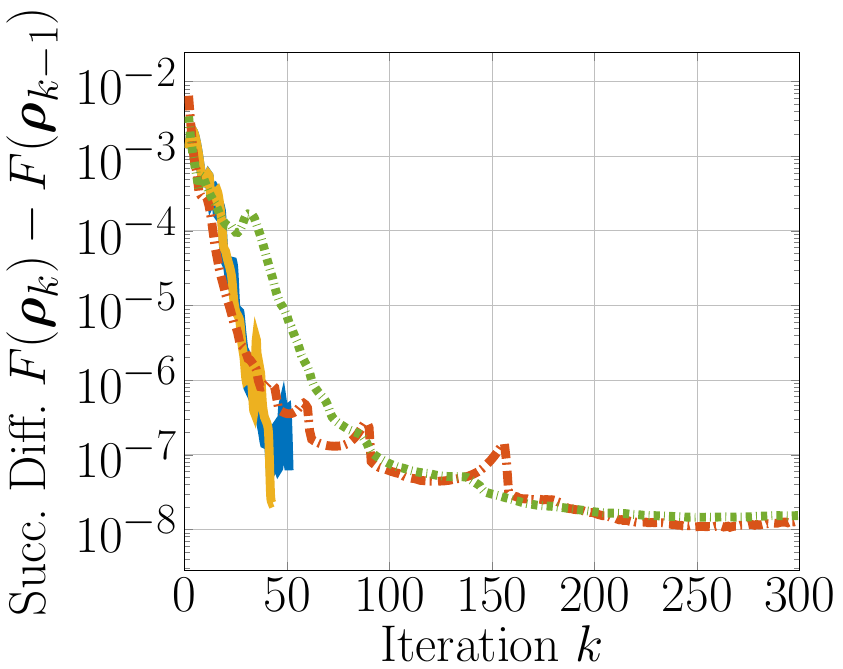}\\[3pt]
  \includegraphics[width=0.375\textwidth]{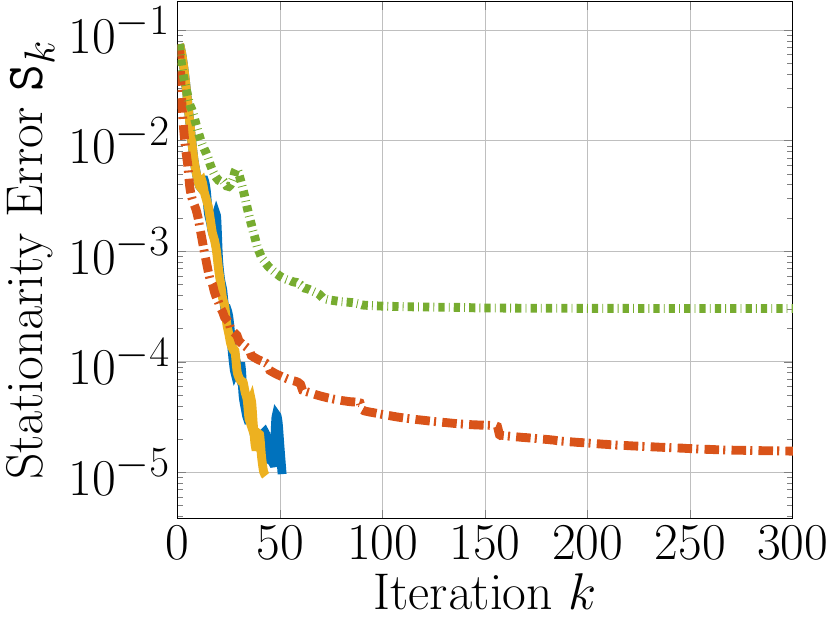}
  ~\;
  \includegraphics[width=0.375\textwidth]{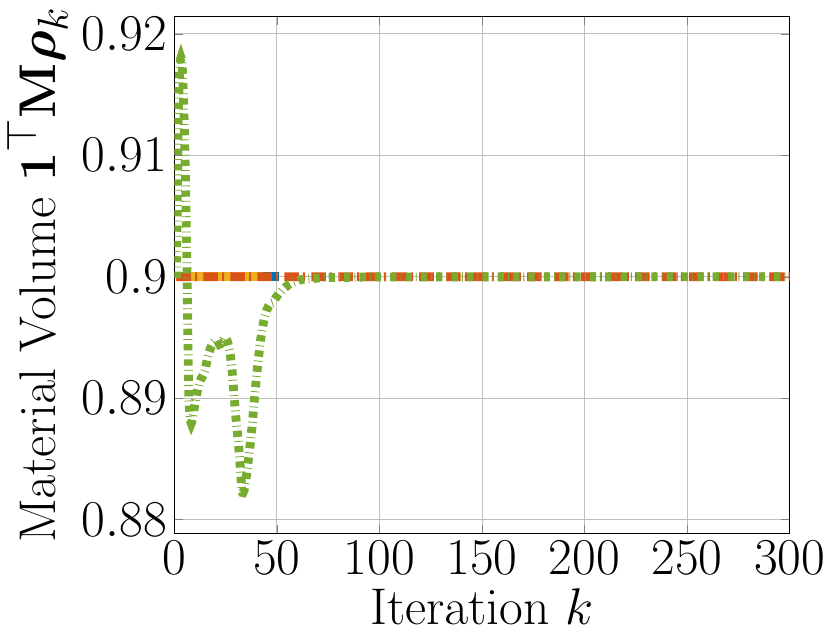}\\
  \caption{Problem 1. Compliance (top left), successive difference of compliance (top right), relative stationarity error (bottom left), and volume (bottom right) for the MBB beam with mesh size $h=1/256$. \label{fig:mbb-compare}}
\end{figure*}
\begin{figure*}
	\centering
  \includegraphics[width=0.315\textwidth]{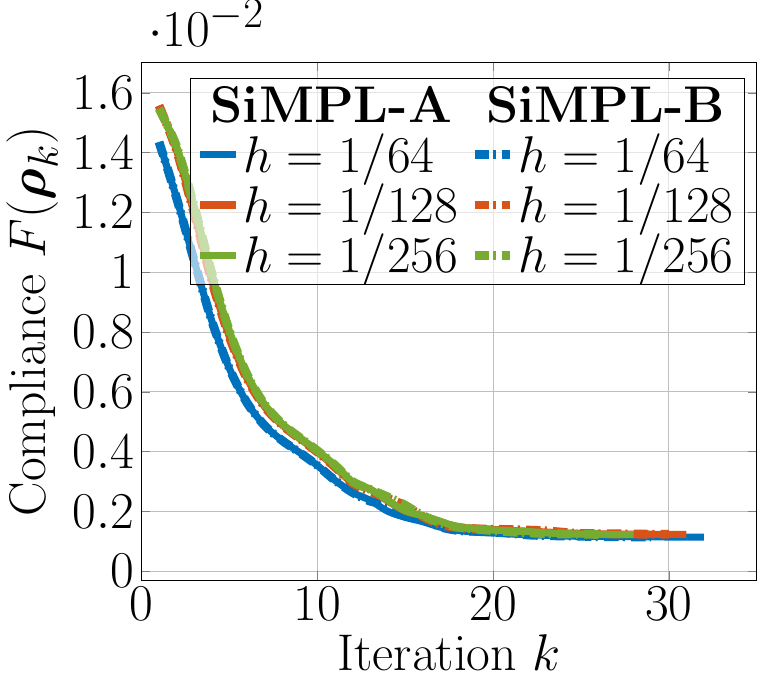}
  \,
  \includegraphics[width=0.33\textwidth]{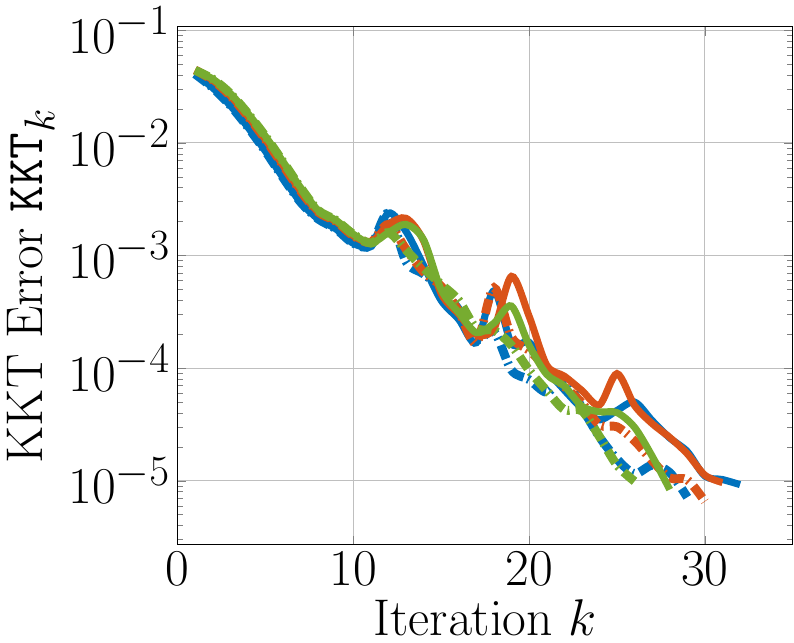}
  \,
  \includegraphics[width=0.318\textwidth]{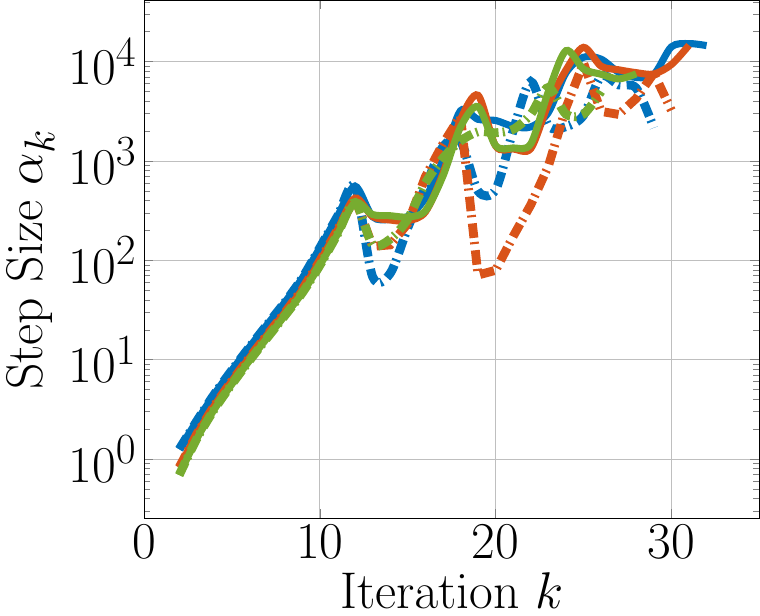}
  \\
  \caption{Problem 1. Compliance (left), relative stationarity error (center), and step size (right) with SiMPL-A and SiMPL-B methods for the MBB beam with various mesh sizes $h=1/64,\;1/128$, and $1/256$.\label{fig:mbb-mesh}}
\end{figure*}

We discretize a $3\times1$ MBB beam into $768\times256$ elements ($h=1/256$) and use a volume fraction of 30\%, i.e.,  $\theta=0.3$.
A horizontal roller supports the bottom right corner of the design domain, and distributed vertical rollers enforce symmetry on the left side of the domain.
An external force $\mathbf{f}$ is applied at the top left corner $\mathbf{c}=(0.0,1.0)^\top$:
\begin{equation*}
  \mathbf{f}=\begin{cases}(0,-1)^\top &\text{if }\|\mathbf{x}-\mathbf{c}\|_{\ell^2},\leq 0.05\\\mathbf{0}&\text{otherwise.}\end{cases}
\end{equation*}

\begin{table}[ht]
  \setlength{\tabcolsep}{4.5pt}
  \renewcommand{\arraystretch}{1.2}
  \centering
  \begin{tabular}{|c||c|c|c|c|c|}\hline
            & $F(\bm{\rho}_{{\rm final}})$ & $\texttt{S}_k$       & Its. & Evals. \\\hline\hline
    SiMPL-A & $1.2078\times 10^{-3}$       & $9.62\times 10^{-6}$ & 50   & 56     \\\hline
    SiMPL-B & $1.2079\times 10^{-3}$       & $9.30\times 10^{-6}$ & 46   & 58     \\\hline
    OC      & $1.2234\times 10^{-3}$       & $1.42\times 10^{-5}$ & 300  & 300    \\\hline
    MMA     & $1.2129\times 10^{-3}$       & $3.01\times 10^{-4}$ & 300  & 300    \\\hline
  \end{tabular}
  \caption{Problem 1. The number of cumulative iterations and objective function evaluations for each method.}
  \label{tab:mbb}
\end{table}
\begin{table*}[ht]
  \setlength{\tabcolsep}{4.5pt}
  \renewcommand{\arraystretch}{1.2}
  \centering
  \begin{tabular}{ccccc}
    \toprule%
    & \multicolumn{4}{@{}c@{}}{$F(\rho_{30}) \textrm{ (volume)}$}  \\\cmidrule{2-5}%
    Mesh size & SiMPL-A& SiMPL-B& OC& MMA\\
    \midrule
    $1/64$ & $1.0322\times 10^{-3}\;(0.90)$ & $1.0147\times 10^{-3}\;(0.90)$ & $1.0665\times 10^{-3}\;(0.90)$ & $1.2407\times 10^{-3}\;(0.89)$\\
    $1/128$ & $1.0941\times 10^{-3}\;(0.90)$ & $1.0831\times 10^{-3}\;(0.90)$ & $1.1415\times 10^{-3}\;(0.90)$ & $1.4883\times 10^{-3}\;(0.88)$\\
    $1/256$ & $1.0566\times 10^{-3}\;(0.90)$ & $1.0801\times 10^{-3}\;(0.90)$ & $1.1336\times 10^{-3}\;(0.90)$ & $2.3581\times 10^{-3}\;(0.88)$\\
    $1/512$  & $1.0989\times 10^{-3}\;(0.90)$ & $1.0972\times 10^{-3}\;(0.90)$ & $1.1397\times 10^{-3}\;(0.90)$ & $2.7049\times 10^{-3}\;(1.01)$ \\
    \botrule
  \end{tabular}
  \caption{Problem 1. Computed objective function values and material volumes (in parentheses) at iteration 30 for different mesh sizes. SiMPL and OC show nearly mesh-independent behavior, while MMA produces mesh-dependent objective function values and volumes.}
  \label{tab:mbb2}
\end{table*}%
\begin{figure*}
	\centering
  \begin{tabular}{cc}
    \includegraphics[width=0.35\textwidth]{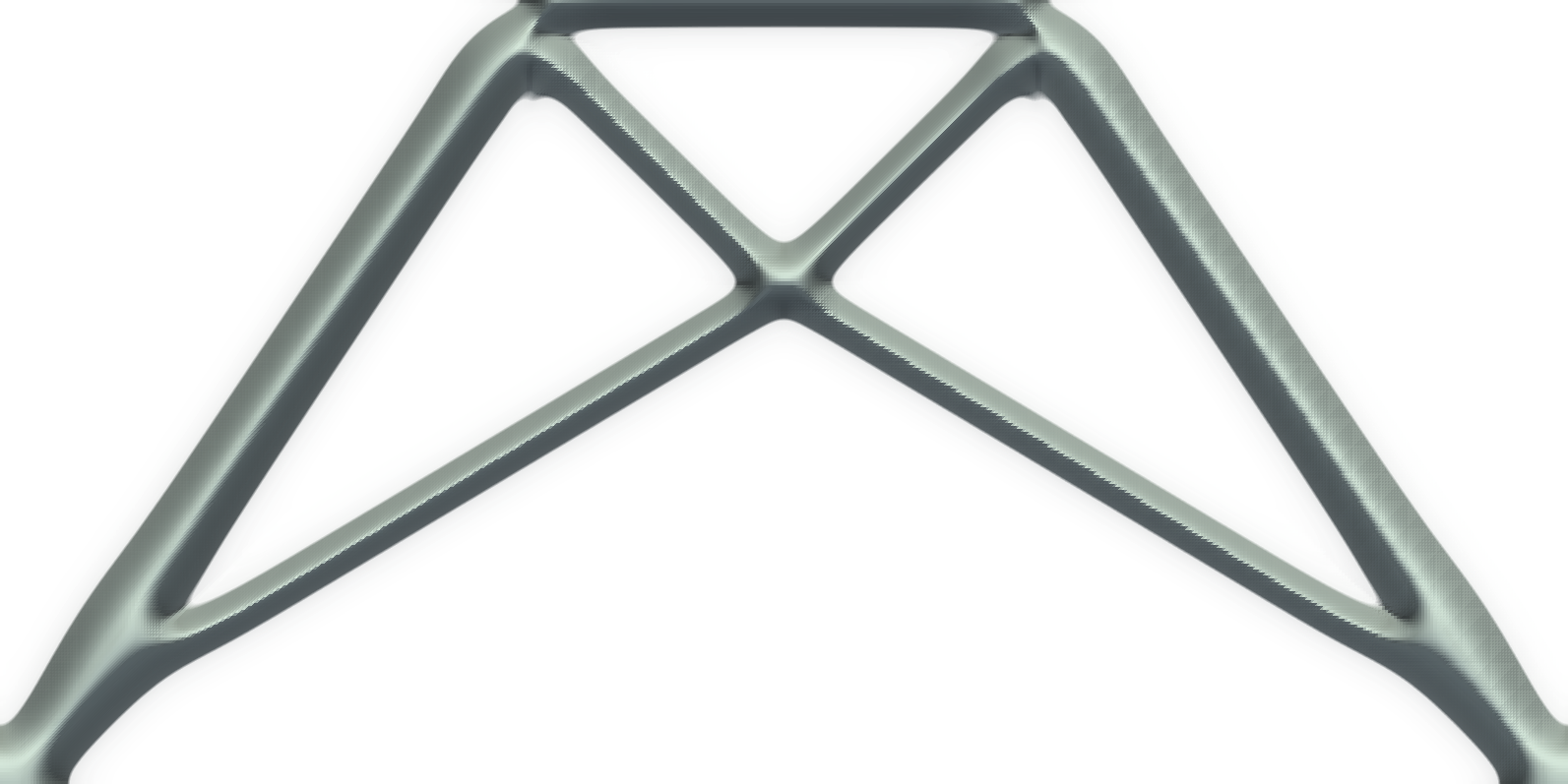}&
    \qquad~
    \includegraphics[width=0.35\textwidth]{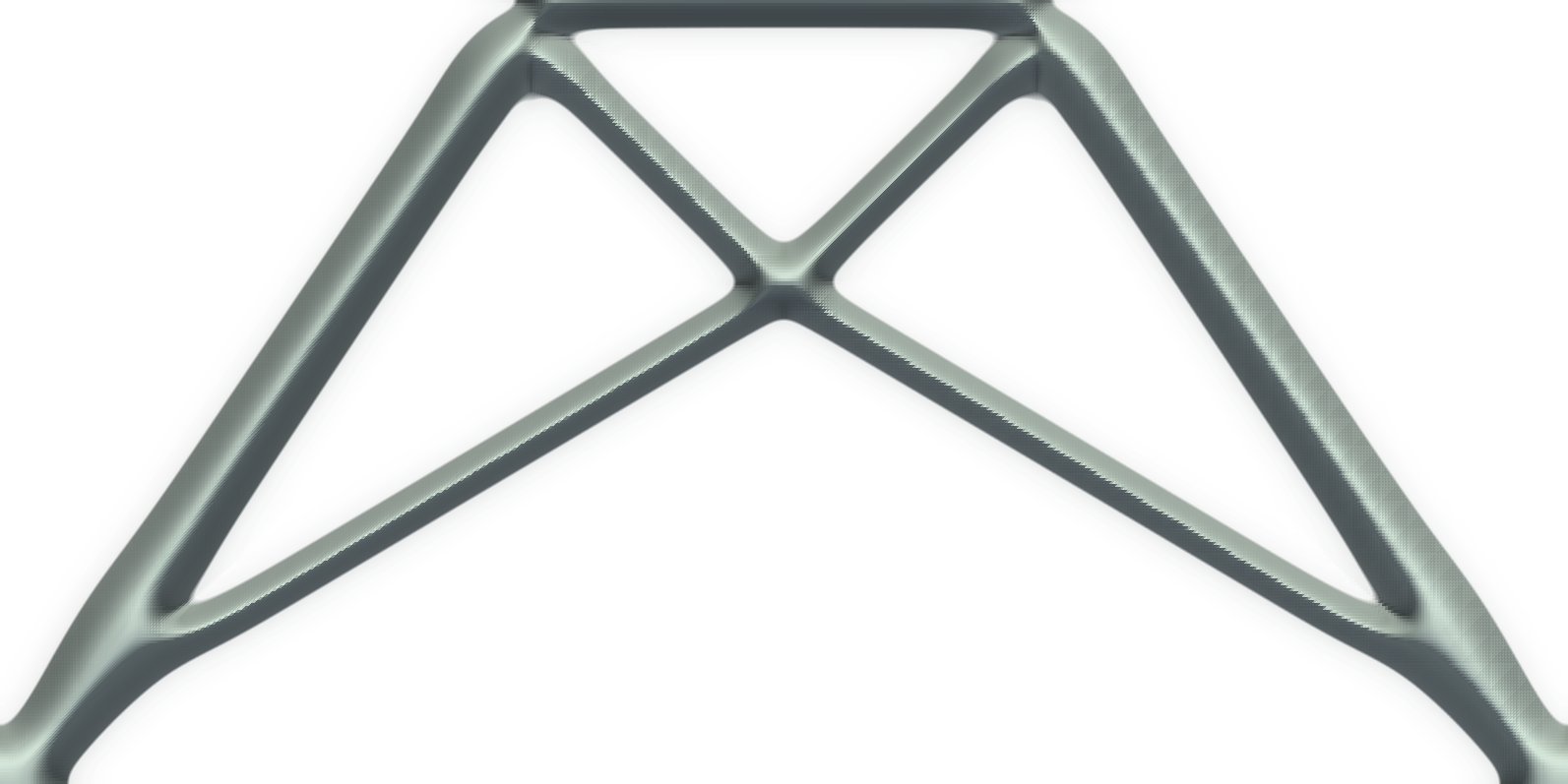}\\
    \textbf{(a)} SiMPL-A&\textbf{(b)} SiMPL-B\\[1em]
  \end{tabular}
  \caption{Problem 3. Optimized frame designs with SiMPL-A (27 iterations) and SiMPL-B (25 iterations) for this multiple load problem. The final objective function values are $2.8092\times 10^{-3}$ and $2.8185\times10^{-3}$, respectively.}
  \label{fig:frame}

  \includegraphics[width=0.30\textwidth]{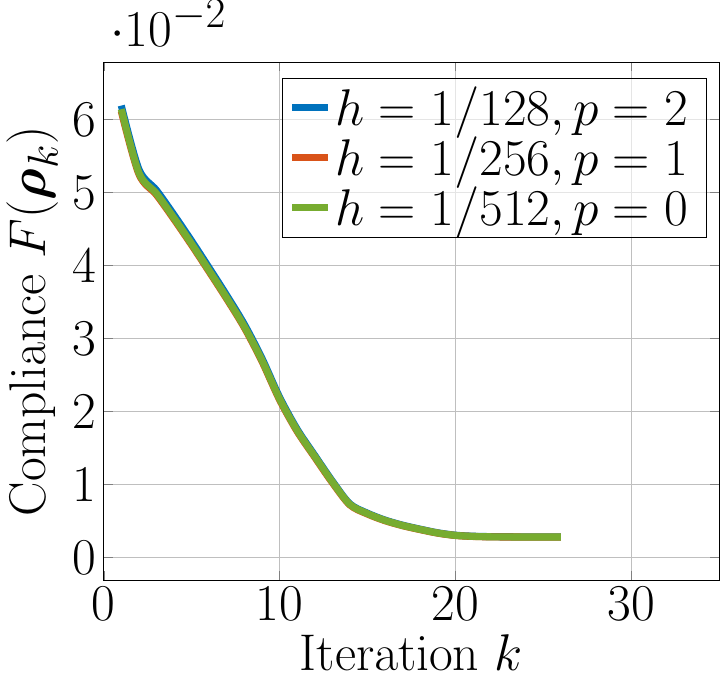}
  \,
  \includegraphics[width=0.333\textwidth]{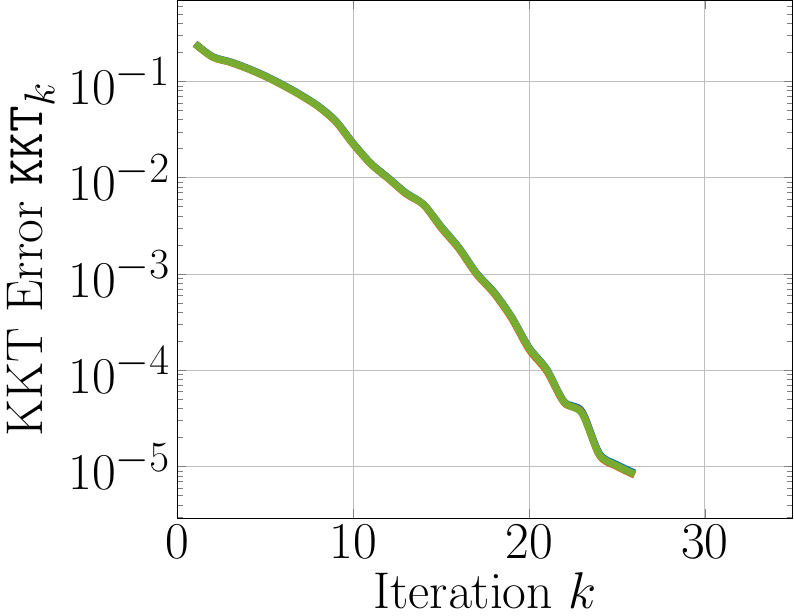}
  \,
  \includegraphics[width=0.33\textwidth]{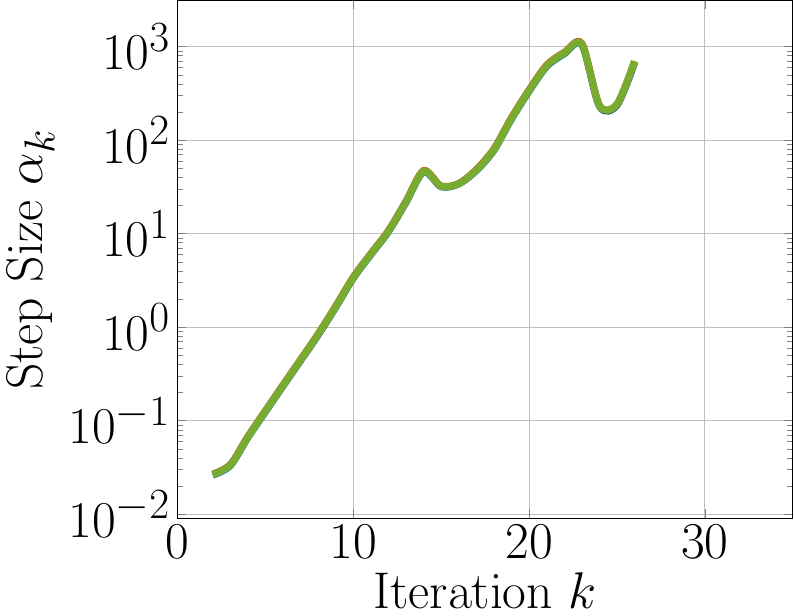}
  \\
  \caption{Problem 3. Compliance (left), relative stationarity error (center), and step size (right) with the SiMPL-B method across various mesh sizes $h=1/128,\;1/256$, and $1/512$ and polynomial degrees $p=0, 1, 2$. The total number of backtracking steps were 3 for both SiMPL-A and SiMPL-B for each mesh resolution.\label{fig:multiload-plot}}
\end{figure*}
\paragraph*{A common stopping criterion}
Our first aim is to compare SiMPL to OC and MMA.
Since the formula for $\bm{\lambda}_k$ in~\eqref{eq:ApproximateLM} is not applicable to the OC and MMA optimization algorithms, we propose a stopping criterion specific to this example.
To this end, we note that the standard first-order optimality conditions for constrained optimization require the negative gradient of the objective function to be a linear combination of the gradients of the active constraint functions \cite[Chap 12]{nocedal1999numerical}.
The coefficients of this linear combination are nothing more than Lagrange multipliers, which satisfy complementarity conditions with the constraints.
Therefore, we propose to compare SiMPL, OC, and MMA by checking the (equivalent) stationary error condition
\begin{equation}
  \label{eq:KKT-L2}
  \mathtt{S}_k
  =
  (\mathbf{s}_k^\top\mathbf{M}\mathbf{s}_k)^{1/2}
  \leq
  10^{-5}
  \,,
\end{equation}
where $\mathbf{s}_k := \bm{\rho_k} - \mathcal{P}(\bm{\rho_k} - \mathbf{g}_k)$.

\paragraph*{Comparison to OC and MMA}
The optimized MBB beam designs for SiMPL-A/B, OC, and MMA are shown in \Cref{fig:mbb}, with convergence histories in \Cref{fig:mbb-compare}.  The resulting topologies are qualitatively similar. A record of the number of iterations and objective function evaluations is given in \Cref{tab:mbb}. SiMPL-A returned a design satisfying~\eqref{eq:KKT-L2} after 50 iterations. On the other hand, SiMPL-B required only 46 iterations but ended up performing more objective function evaluations because it also took more backtracking steps.
Both OC and MMA appear to require more than 300 iterations to reach the prescribed stationarity error \eqref{eq:KKT-L2}. Even at iteration 300, the compliance values for OC and MMA are larger than with SiMPL-A or SiMPL-B. 
In \Cref{tab:mbb2}, we record the computed objective function values and material volumes at iteration 30 for different mesh sizes.
The results show that SiMPL and OC are nearly mesh-independent, while MMA exhibits mesh-dependent behavior, with the accuracy reducing each time the mesh is refined.

\paragraph{Mesh-independence}
In the previous example, we observed that the SiMPL method shows nearly mesh-independent behavior.
Taking the same 2D MBB beam problem, we now investigate and demonstrate the mesh-independent behavior of the SiMPL method.
Here, we utilize the stopping criterion proposed in \eqref{eq:KKT_estimator} and rerun SiMPL-A and SiMPL-B on discretizations with mesh sizes $h=1/64,1/128,1/256$.
\Cref{fig:mbb-mesh} depicts the resulting objective function values, KKT errors, and step sizes throughout the optimization process.
Each plot shows similar behavior indicating mesh-independence of the method.
In Problem 3, we investigate degree-independent behavior of the SiMPL method with a more complex example.

\bigskip
As a second compliance minimization example, we consider optimizing a 3D cantilever beam across different material volumes.

\paragraph*{Problem 2: 3D cantilever beam}

\begin{figure*}[!ht]
  \centering
  \begin{tabular}{cc}
    \includegraphics[width=0.35\textwidth]{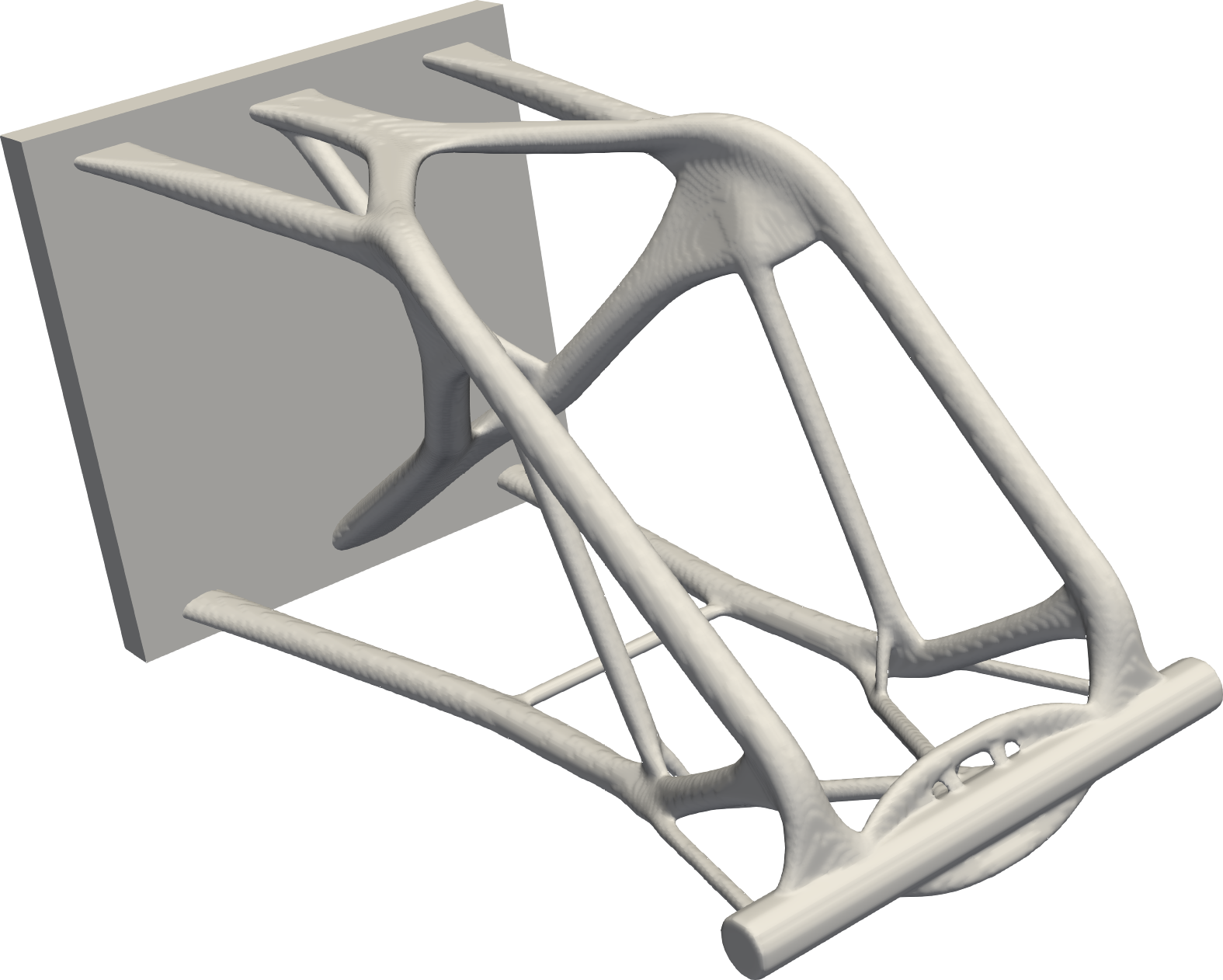}&
    \qquad~
    \includegraphics[width=0.35\textwidth]{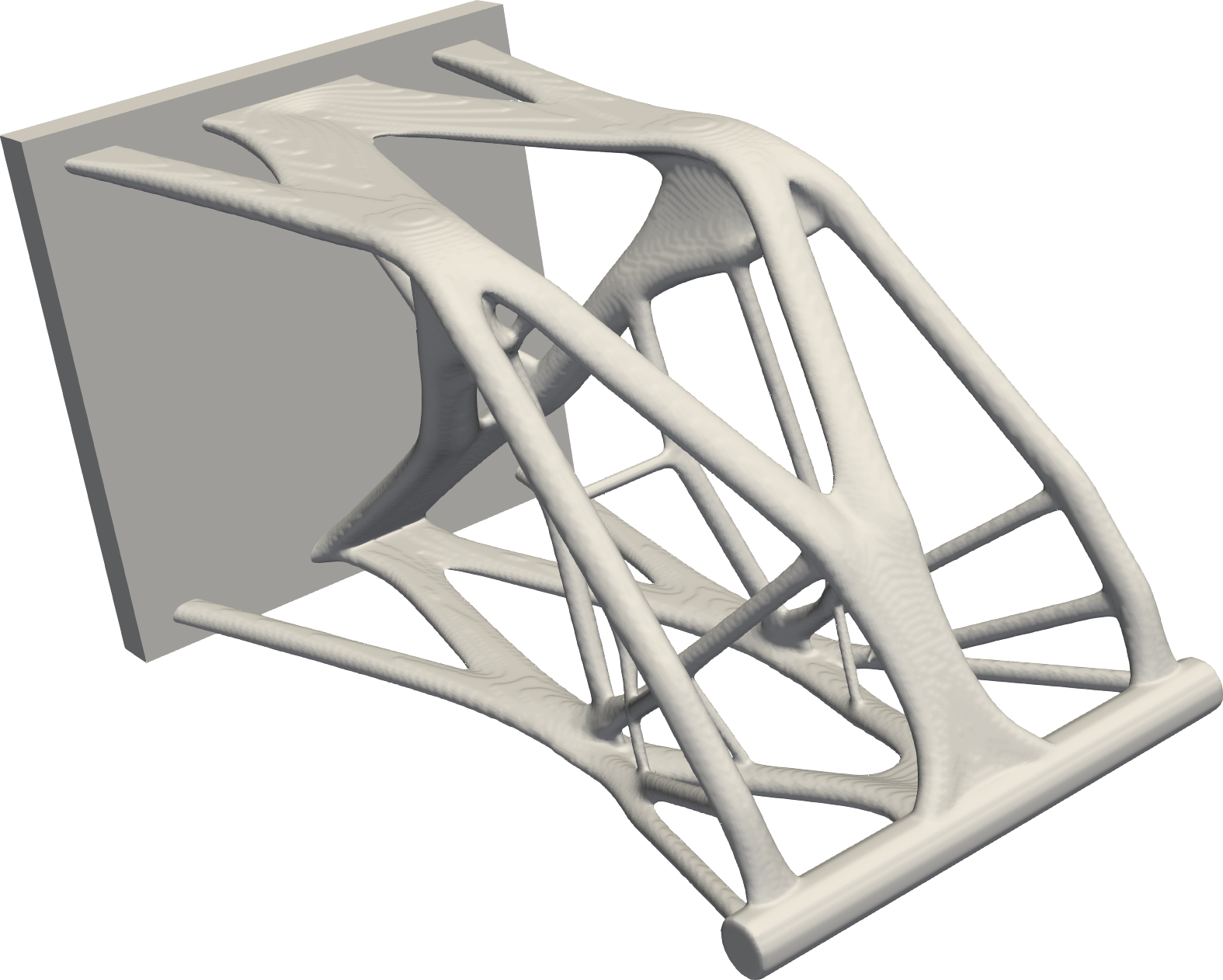}\\
    \textbf{(a)} $\theta=0.07$&\textbf{(b)} $\theta=0.10$\\[1em]
    \includegraphics[width=0.35\textwidth]{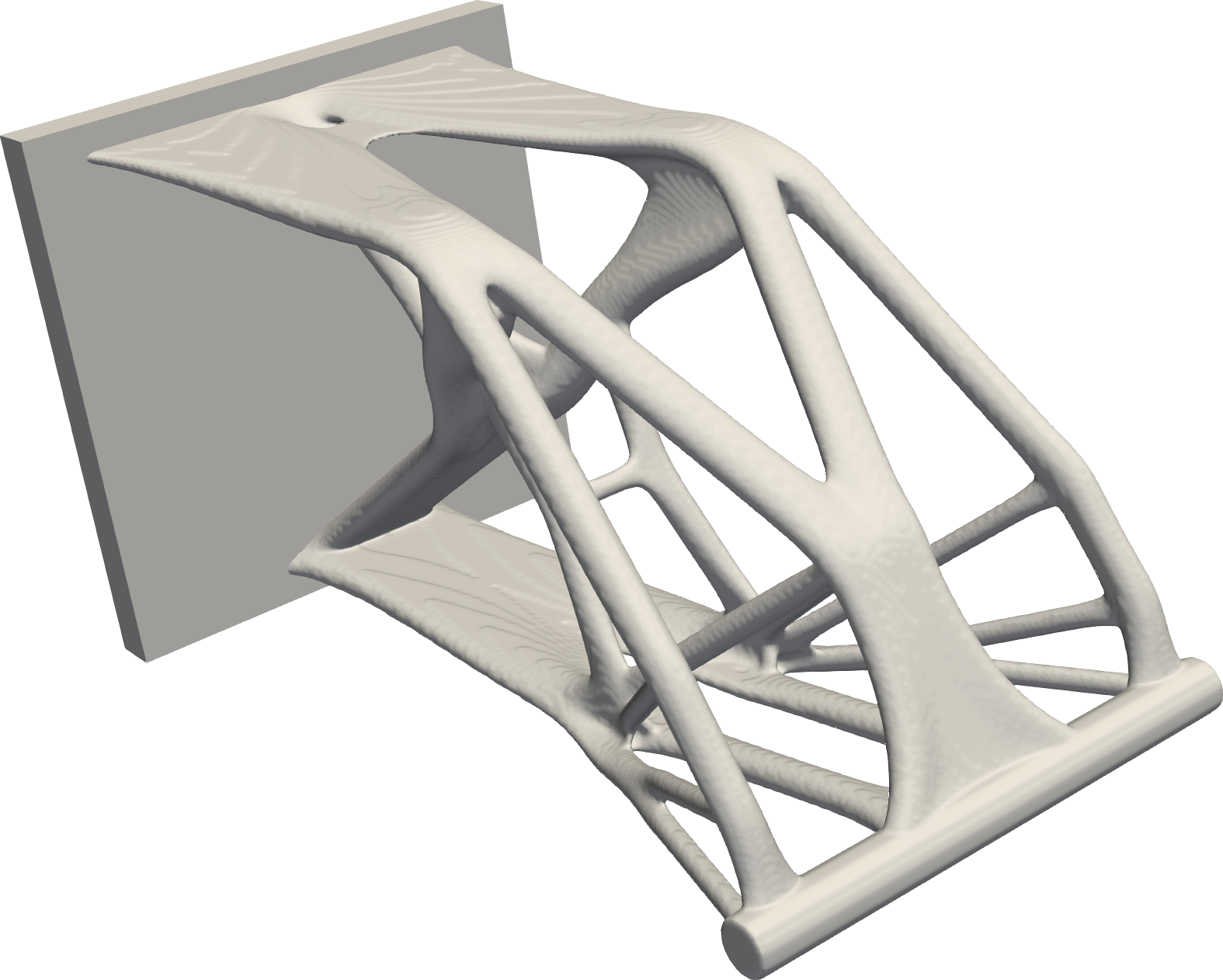}&
    \qquad~
    \includegraphics[width=0.35\textwidth]{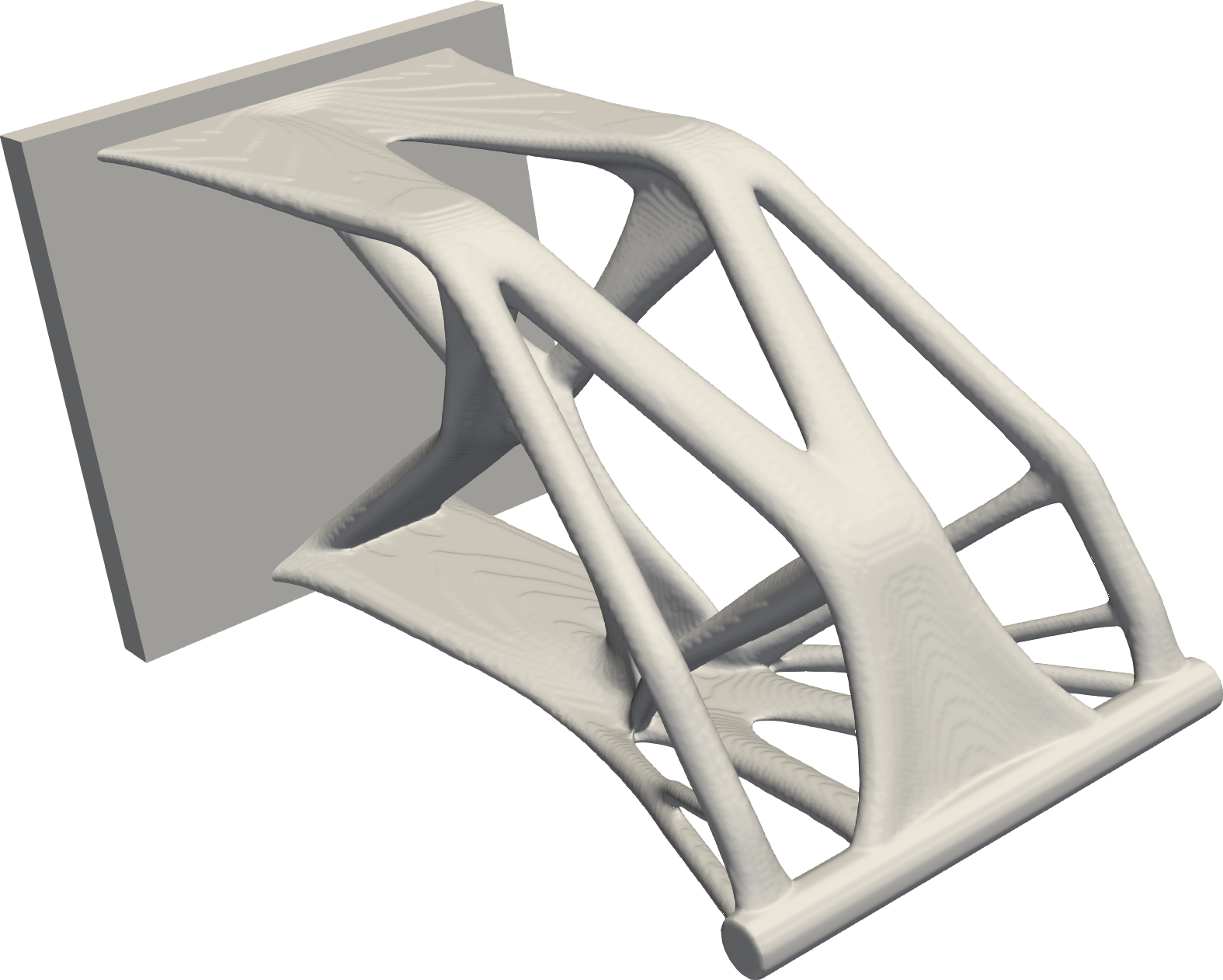}\\
    \textbf{(c)} $\theta=0.12$&\textbf{(d)} $\theta=0.15$\\[1em]
    \includegraphics[width=0.35\textwidth]{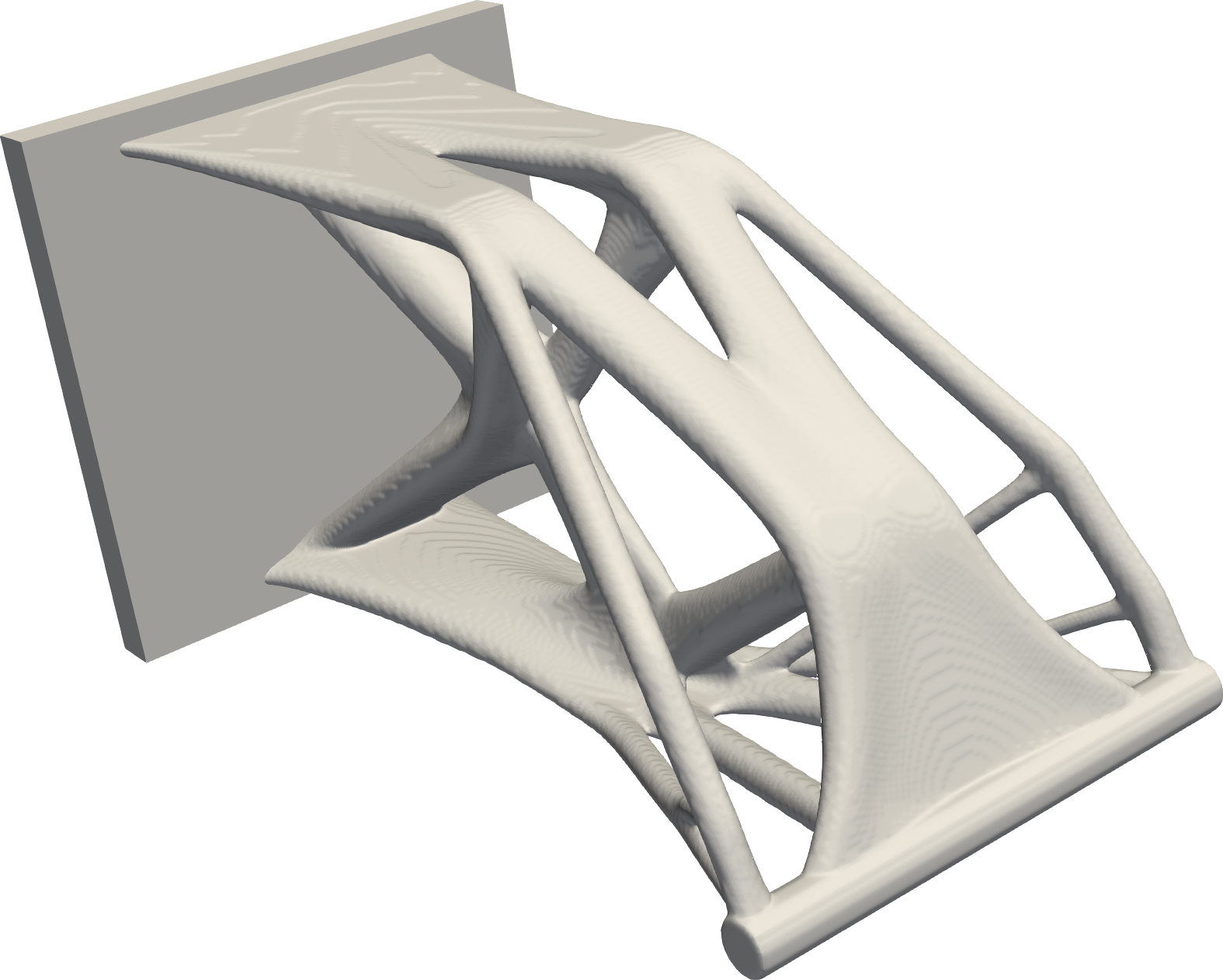}&
    \qquad~
    \includegraphics[width=0.35\textwidth]{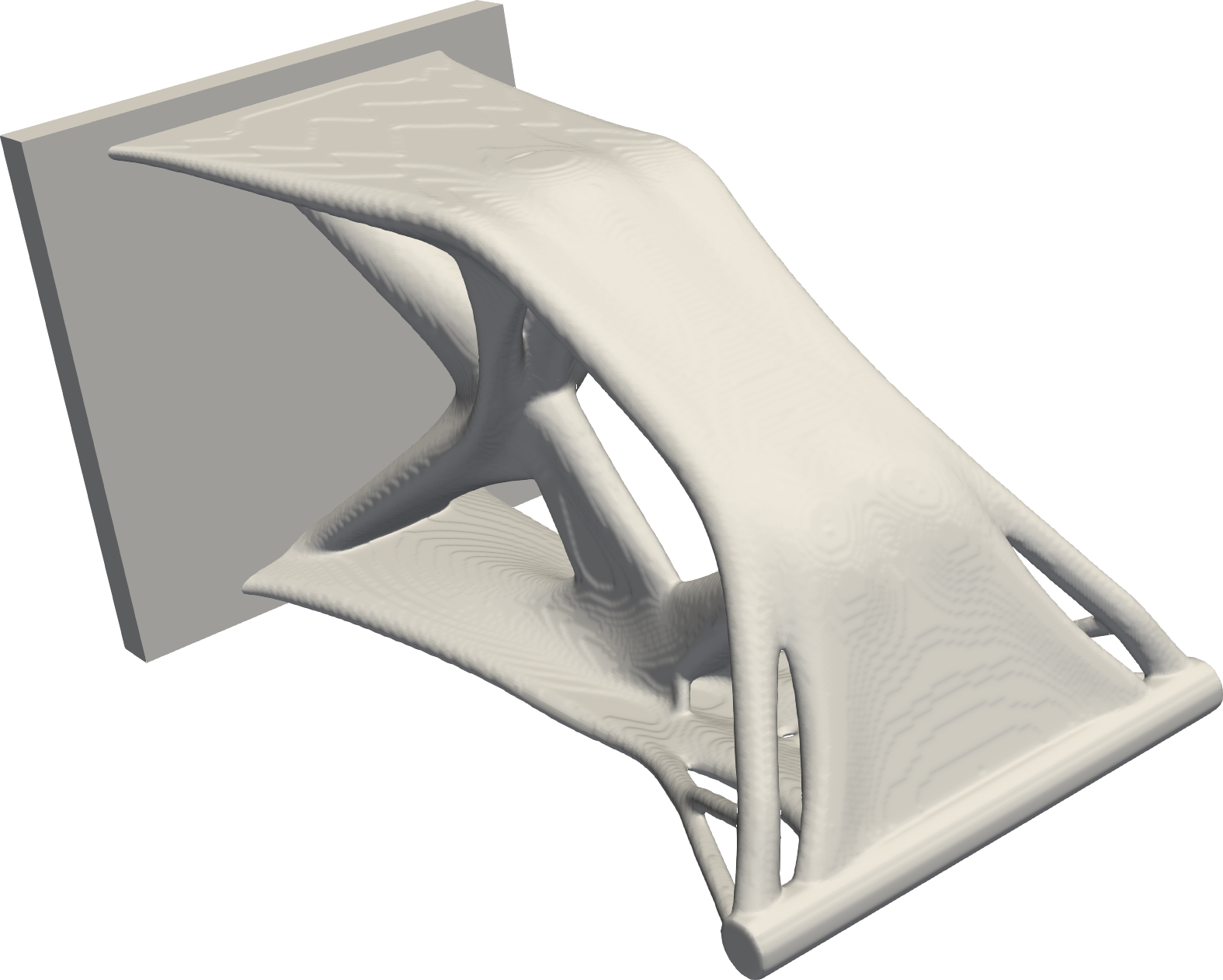}\\
    \textbf{(e)} $\theta=0.17$&\textbf{(f)} $\theta=0.20$
  \end{tabular}
  \caption{Problem 2. The iso-surfaces $\{ x \mid \tilde{\rho}(x) =0.5 \}$ with SiMPL converged in 81, 44, 48, 44, 49, and 42 iterations with $\theta=[0.075,0.10,0.125,0.15,0.175,0.20]$ and objective function values $\left[8.04, 5.71, 3.96, 3.03, 2.51, 2.01\right]\times10^{-2}$, respectively.
    \label{fig:canti3}
  }
\end{figure*}

We discretize $2\times1\times1$ cantilever into $512\times256\times256$ elements. A downward distributed load is applied at passive solid elements $\{(x-1.9)^2+(z-0.1)^2<0.05^2\}$ and the left side of the design domain $\{x=0.0\}$ is fixed in all displacement components.
We then run SiMPL-B with $\mathtt{tol} = 10^{-5}$ for volume fractions $\theta$ varying from $0.075$ to $0.2$.
The optimized designs, shown in \Cref{fig:canti3}, required between 81 ($\theta = 0.075$) and 42 ($\theta = 0.2$) iterations. The wide variation between the required numbers of iterations is due to differences in the complexities of the final designs.
In particular, an optimized design with a volume fraction of 7.5\% cannot form large solid members, and the optimized design contains an interconnected set of beam-like members. Such a topology differs significantly from the one appearing in the early algorithm iterations, leading to the complete removal of several features and subsequent adjustments of the remaining ones, especially around the load region.
Notably, all of the designs with $0.1 \leq \theta \leq 0.2$ converged after a similar number of iterations.

\subsection{Multiple load compliance minimization}
The next application is the compliance minimization problem with multiple external loads.
\begin{align*}
  \min_{\bm{\rho}}\  & \sum_{l=1}^{N_\ell}(\mathbf{f}^l)^\top (\mathbf{u}^l)
  \\
  \text{subject to }~
                     &
                     \mathbf{K}(\tilde{\bm{\rho}})\mathbf{u}^l =\mathbf{f}^l,\quad l=1,...,N_\ell,                                                  \\
                     & (\epsilon^2\mathbf{A}+\tilde{\mathbf{M}})\tilde{\bm{\rho}} =\mathbf{N}\bm{\rho}  \\
                     & \mathbf{0}\leq \bm{\rho} \leq \mathbf{1},                                        \\
                     & \mathbf{1}^\top\mathbf{M}\bm{\rho} \leq\theta|\Omega|
                     \,,
\end{align*}
where $N_\ell$ is the number of external loads, and $\mathbf{f}^l$ and $\mathbf{u}^l$ are the external force and displacement corresponding to the $l$-th load, respectively.

In here, we consider the sum (average) of compliance for external loads.
This problem can be extended to the worst-case compliance minimization problem by considering the maximum compliance among external loads.
However, this extension is beyond the scope of this paper, and we leave it for future work.

\paragraph*{Problem 3: 2D frame}
\normalfont
We seek an optimized frame structure on a $2\times1$ domain with a volume fraction of 20\%.
Both the left and right bottom corners of the domain are pin-supported.
Two external vertical forces ($N_\ell=2$) are applied at the one-third and two-thirds points of the top boundary of the domain:
\begin{equation*}
  \mathbf{f}^l=\begin{cases}(0,-1)^\top &\text{if }\|\mathbf{x}-\mathbf{c}_l\|_{\ell^2}\leq 0.05\\\mathbf{0}&\text{otherwise.}\end{cases}
\end{equation*}
Here, $\mathbf{c}_1=(2/3,0.9)^\top$ and $\mathbf{c}_2=(4/3,0.9)^\top$.
In this example, we compare the convergence histories of SiMPL-B for polynomial orders $p=0,1,2$ in the unfiltered density space.
When $p\geq1$, the filtered density and displacement fields are computed using $p$-th order $C^0$-conforming polynomials.
For $p=0$, these fields are computed using a piecewise-linear continuous space to maintain the conformity of the discrete spaces.
The filter solver, discretized with standard finite element method, does not guarantee the discrete maximum principle, and hence $0\leq \tilde{\rho}_h\leq 1$ may not be satisfied.
To ensure the well-posedness of the elasticity equation, we clipped the discrete filtered density to the range $[0,1]$.
While we have not observed significant impact on the optimization process, clipping can be avoided by using a finite element method that satisfies the discrete maximum principle, e.g., \cite{keith2023proximal, fu2024}.
The optimized design for each $p$ and $h$ is depicted in \Cref{fig:frame}.
Under the same stopping criterion $\mathtt{KKT}_k\leq 10^{-5}$, SiMPL-B converged in 25 iterations with 28 objectives evaluation for all $p$ values.
The plots in \Cref{fig:multiload-plot} illustrate the behavior of the objective function, stationarity error, and step size for different polynomial degrees in the finite element space.
Notably, all polynomial degrees exhibit nearly identical trends, suggesting that the SiMPL method is degree-independent.
\subsection{Self-weight compliance minimization}
\begin{figure*}[t]
  \centering
  \includegraphics[width=0.9\textwidth]{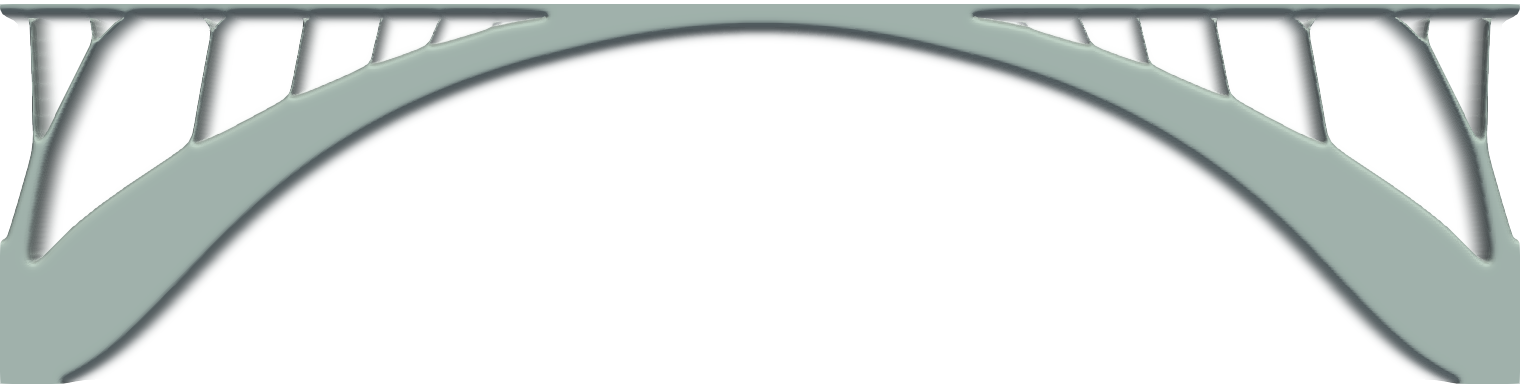}
  \caption{Problem 4. An optimized bridge design with SiMPL-B converged in 81 iterations with $F(\rho_{\rm final})=779.99$. The design is optimized in the half domain $(0,2)\times(0,1)$ with $1024\times 512$ elements and reflected about $x=0$ for visualization.}
  \label{fig:bridge}
\end{figure*}

The next considered application is the problem of self-weight compliance minimization:

\begin{equation*}
  \begin{aligned}
    \min_{\bm{\rho}}\  & \mathbf{f}^\top\mathbf{u} + \mathbf{g}(\tilde{\bm{\rho}})^\top\mathbf{u}                        \\
    \text{subject to }
                         & \mathbf{K}(\tilde{\bm{\rho}})\mathbf{u}             =\mathbf{f} + \mathbf{g}(\tilde{\bm{\rho}}) \\
                         & (\epsilon^2\mathbf{A}+\tilde{\mathbf{M}})\tilde{\bm{\rho}}  = \mathbf{N}\bm{\rho}                         \\
                         & \mathbf{0}\leq \bm{\rho}                            \leq \mathbf{1},                                      \\
                         & \mathbf{1}^\top\mathbf{M}\bm{\rho}                     \leq\theta|\Omega|
                         \,,
  \end{aligned}
\end{equation*}
where $\mathbf{f}$ is an external force and $\mathbf{g}(\tilde{\bm{\rho}})$ is a downward internal force with magnitude $9.81(\tilde{\bm{\rho}})_i$ at each element $i$.

\paragraph*{Problem 4: Self-weighted bridge}
\normalfont
We seek an optimized bridge on a $2\times1$ domain partitioned into $1024\times 512$ elements with roller boundary conditions on the left-hand side of the domain to enforce a symmetric design.
In addition we choose the volume fraction $\theta=0.7$.
The bridge is pin-supported at the bottom-right corner of the domain and a narrow band of passive elements are used at the top of the domain, $\{(x,y) \mid y\geq 1 - 2^{-5}\}$, to apply a downward force $\mathbf{f}$ with magnitude 40.
The final design, obtained with SiMPL-B and achieving the relative stopping criterion $\texttt{KKT}_k\leq 10^{-5}\texttt{KKT}_0$ after 81 iterations with 52 backtracking steps, is depicted in \Cref{fig:bridge}.
In this case, we found the volume constraint~\eqref{eq:volume-constraint} was inactive with $\mathbf{1}^\top\mathbf{M}\bm{\rho}_{\rm final} = 0.5415 < 0.7|\Omega|$.

\begin{figure*}
    \centering
    \includegraphics[width=0.32\textwidth]{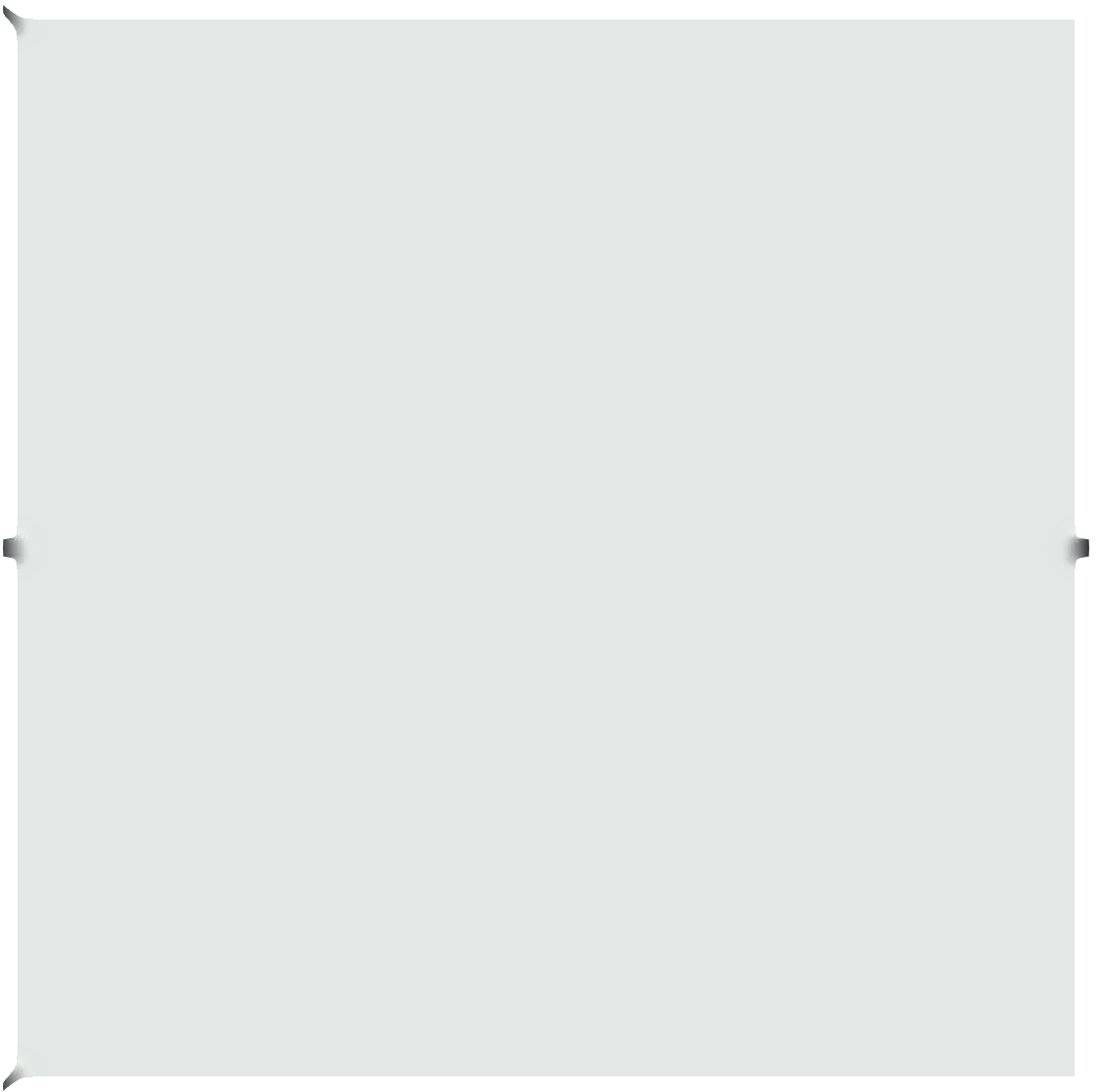}
    \includegraphics[width=0.32\textwidth]{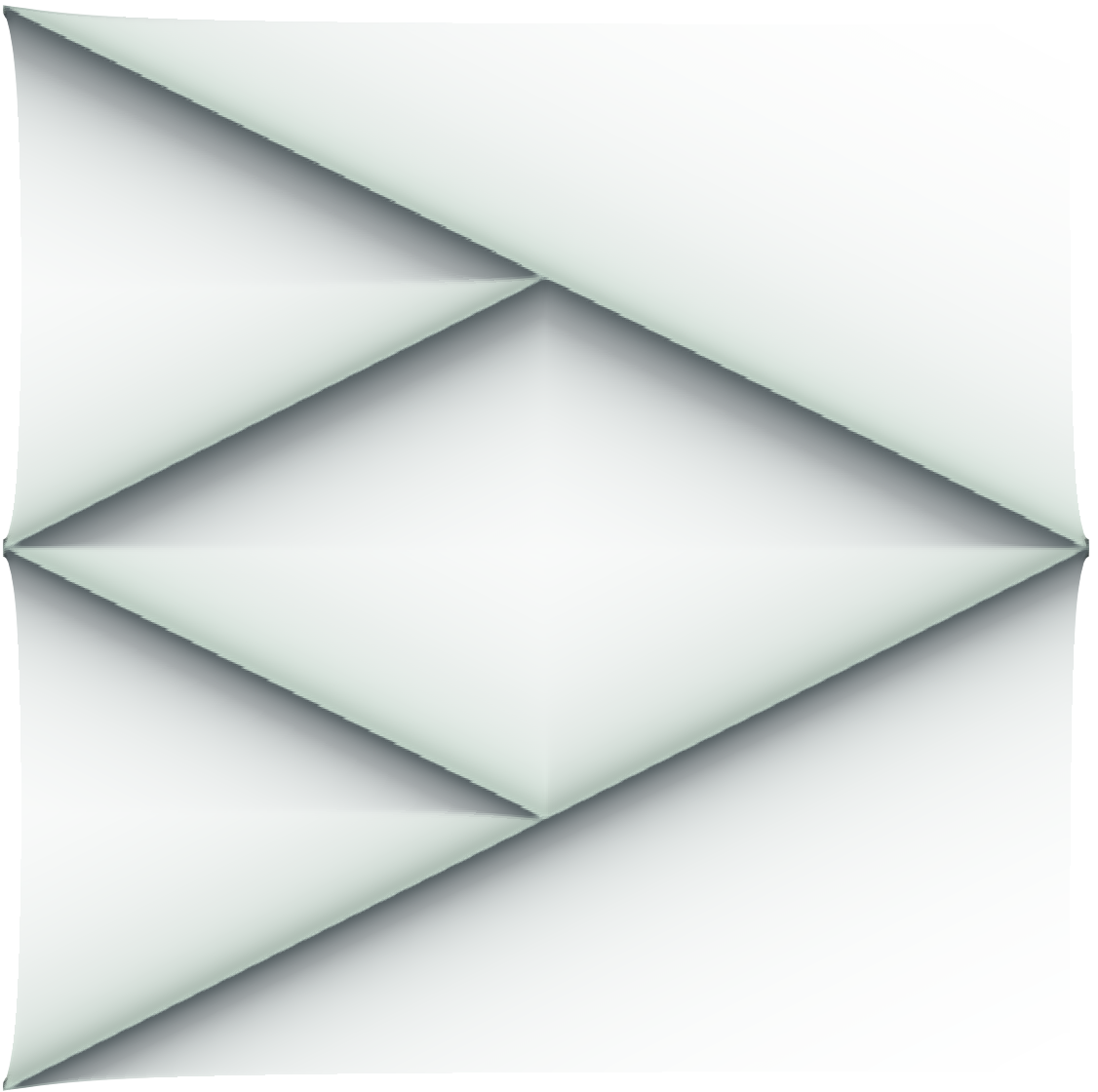}
    \includegraphics[width=0.32\textwidth]{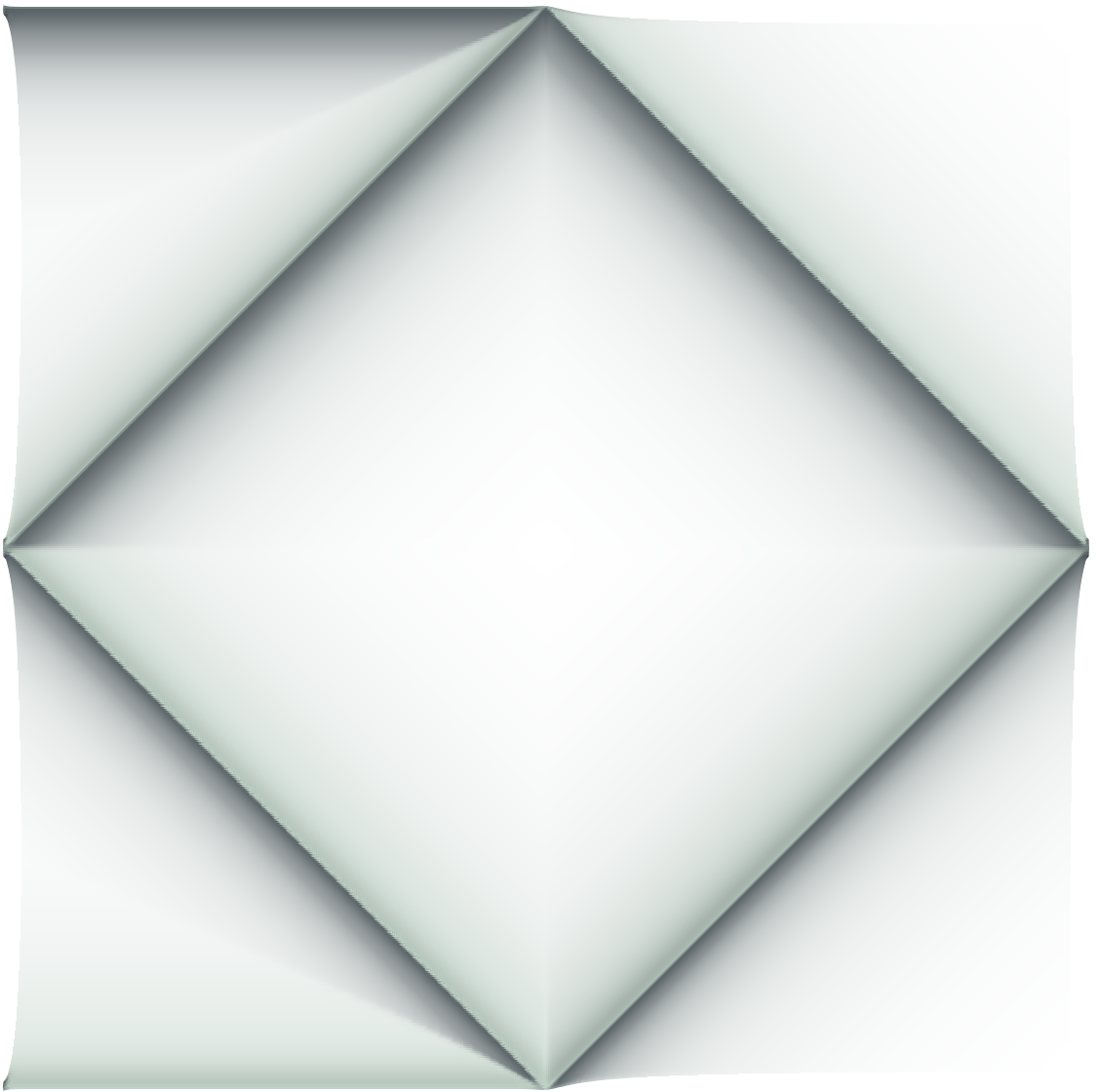}
    \includegraphics[width=0.32\textwidth]{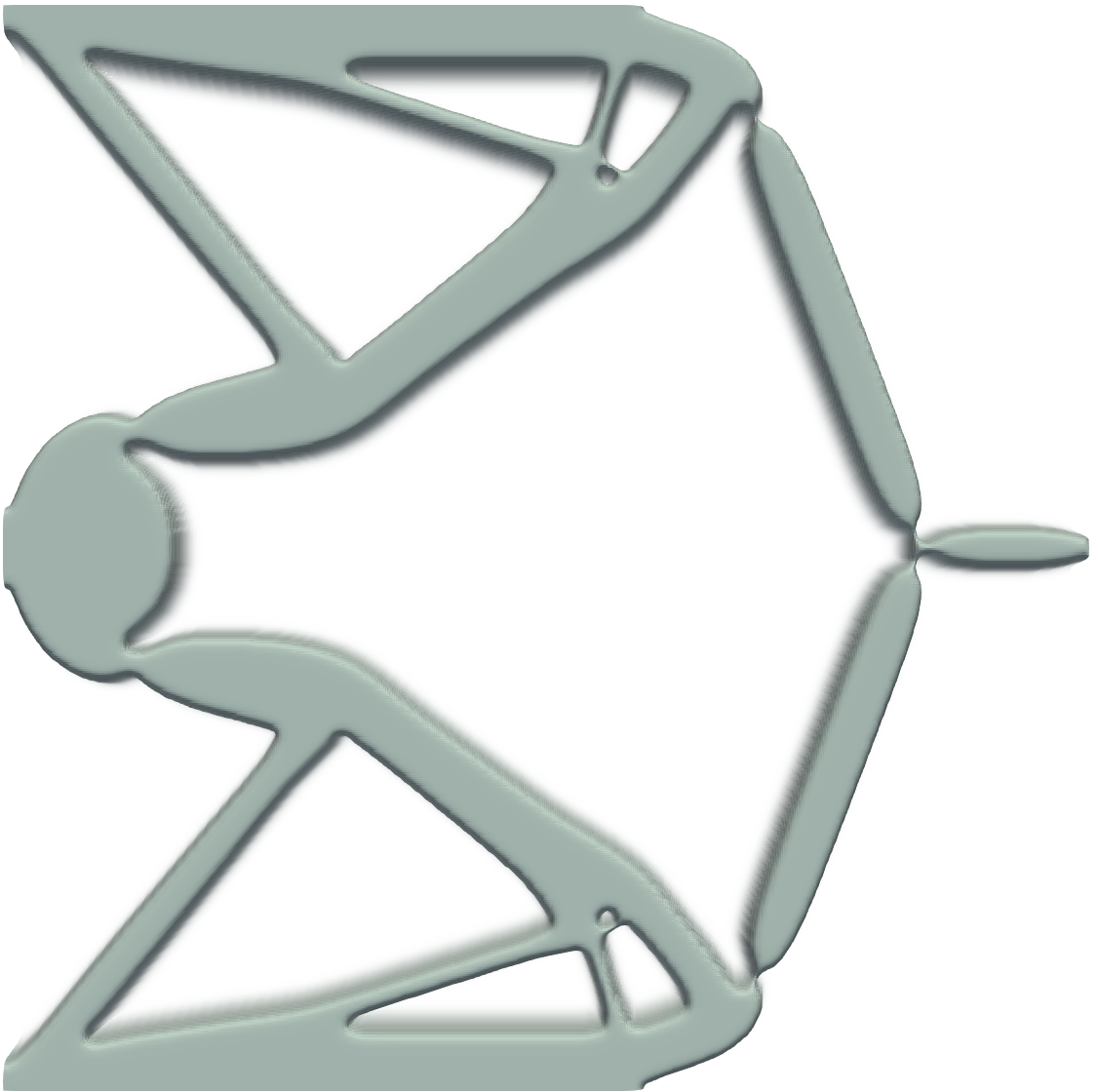}
    \includegraphics[width=0.32\textwidth]{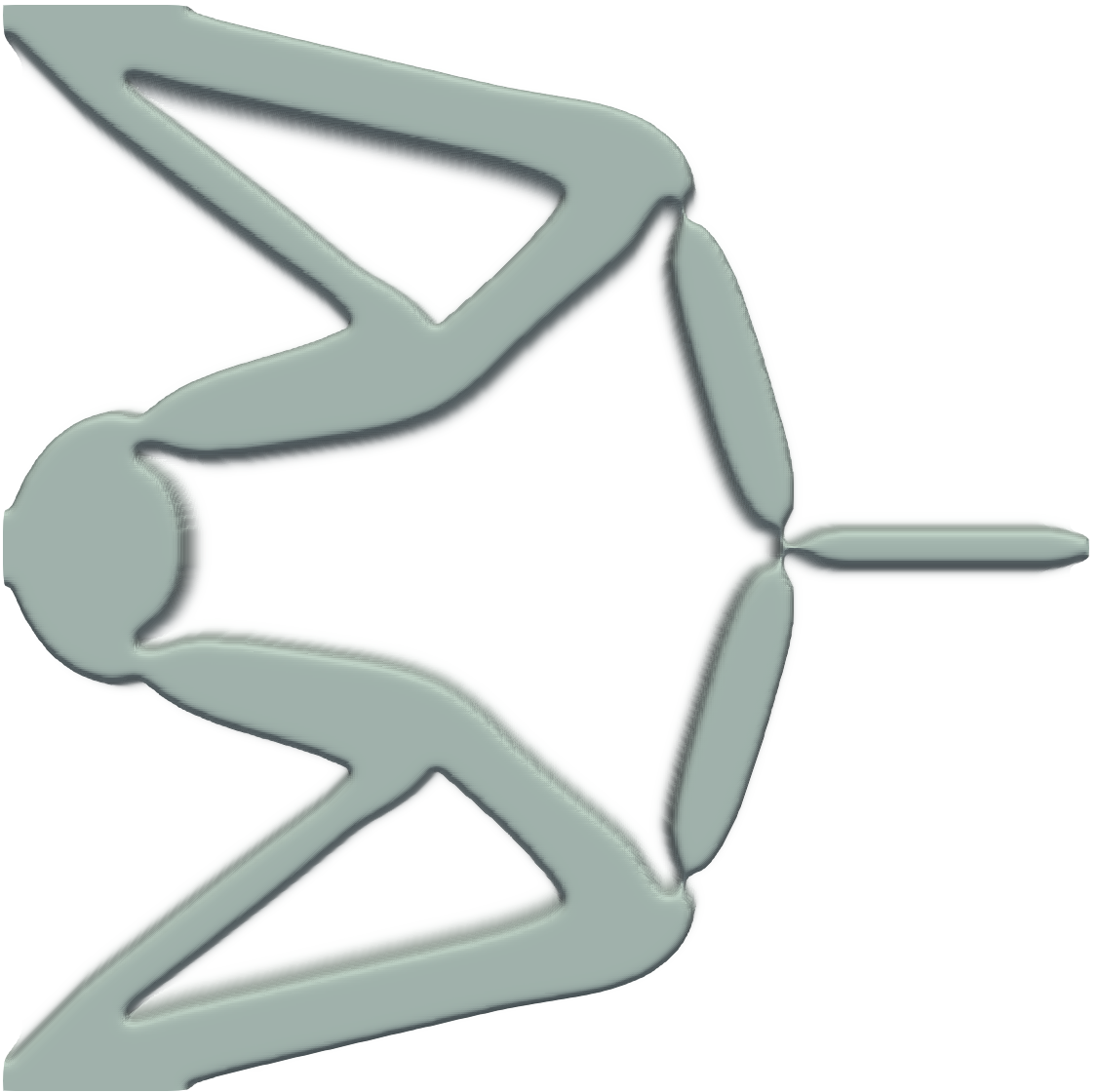}
    \includegraphics[width=0.32\textwidth]{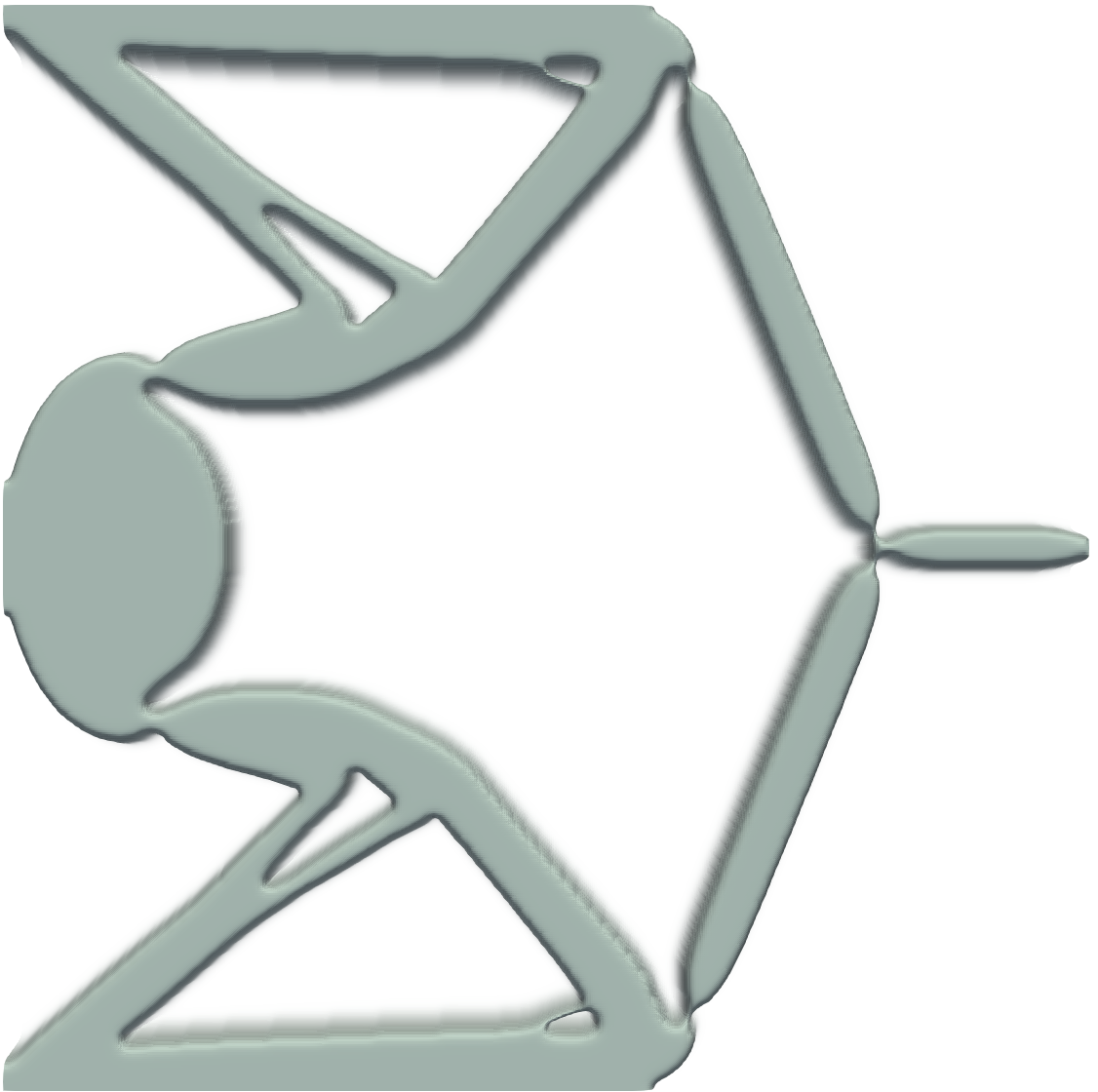}
    \caption{
      Problem 5.
      Initial designs (top row) and optimized designs (bottom row).
      The optimized designs, obtained after 197 (left), 124 (middle) and 139 (right) iterations, have objective function values of $F(\rho_{\rm final})=-0.4598,\;-0.4736$, and $-0.4623$, respectively.
    }
    \label{fig:inverter}
  \end{figure*}

\subsection{Compliant mechanism}
The final application is the compliant mechanism design problem, where the objective is to maximize the displacement resulting from a given input force.
We used the spring and load model \cite{Bendsoe2004}:
\begin{equation*}
  \begin{aligned}
    \min_{\bm{\rho}}\  &-\frac{k_{\rm out}}{L}\mathbf{r}_{\rm out}^\top \mathbf{u}                                                    \\
    \text{subject to }
                       &\mathbf{K}(\tilde{\bm{\rho}})\mathbf{u}            =\frac{k_{\rm in}}{L}\mathbf{d}_{\rm in}             \\
                  &(\epsilon^2\mathbf{A}+\tilde{\mathbf{M}})\tilde{\bm{\rho}} =\mathbf{N}\bm{\rho} \\
                  &\mathbf{0}\leq \bm{\rho}                                         \leq \mathbf{1},                       \\
                  &\mathbf{1}^\top\mathbf{M}\bm{\rho}                    \leq\theta|\Omega|
                  \,.
  \end{aligned}
\end{equation*}
Here, $\mathbf{r}_{\rm out}$ and $\mathbf{d}_{\rm in}$ are vectors corresponding to surface integrals over the input/output ports in the input/output force directions, respectively.
The input and output ports are segments of length $L$.
We set the spring coefficients to be $k_{\rm in}=1.0$ and $k_{\rm out}=0.0005$ and the volume fraction to be $\theta = 0.3$.

\paragraph*{Problem 5: Force inverter}
  We replicate the force inverter problem in \cite{Wang2011} on a square domain partitioned into $512\times 512$ elements.
  The objective of this problem is to maximize the displacement at the output port on the middle right-hand side of the mechanism in the direction opposite to the input force.
  Both the input and output ports are modeled by eight-element-long segments, each with length $L=1/64$, on the middle of the left-hand and right-hand boundaries of the domain, respectively.
  The left-hand corners of the domain are pin-supported, and an inward force is applied parallel to the $x$-axis at the input port.
  We exploited the symmetry of the problem to reduce the computational cost.
  Note that the optimal design depends on $\bm{\rho}_0$ because the SiMPL method can only find a locally-optimal design.
  Thus, in this experiment, we compared SiMPL-B's performance across three different initial design choices.
  In addition to the constant initial design $\bm{\rho}_0 = \theta \mathbf{1}$, two non-uniform initial designs were considered by increasing the density along lines connecting the input/output ports and an intermediary point in the domain.
  These initial designs and the corresponding optimized designs are depicted in \Cref{fig:inverter}.
  Here, we used the complementarity vector $\bm{\eta}_k$ in \eqref{eq:kkt-thomas} and set stopped the algorithm once $\texttt{KKT}_k\leq 5\times 10^{-5}$.
  The number of backtracking steps were $47, 51$ and $51$ for the three designs, respectively.
  All designs exhibit \textit{de-facto hinges}, which can be avoided by adding more design constraints \cite{Wang2011,Lazarov2016}.
\section{Concluding remarks}\label{sec:discussion}
In this work, we introduce the SiMPL method for density-based topology optimization.
The derivation emphasizes the most common discretization choice used in applications and highlights important practical features such as selecting step sizes and stopping criteria.
To this end, we suggest two backtracking line search strategies, Armijo and Bregman, and recommend a stopping criterion derived from the KKT optimality conditions.
Combining these features with an initial step size guess coming from a local estimate of the relative smoothness of the reduced objective function, the SiMPL method yields lower complexity and faster convergence than popularized methods such as OC and MMA.
The most unique feature of the SiMPL method is its use of a latent variable $\bm{\psi}_k$, which represents the (bounded) design density $\mathbf{0} \leq \bm{\rho}_k \leq \mathbf{1}$ in an unbounded space.
The resulting relationship $\bm{\rho}_k = \sigma(\bm{\psi}_k)$, where $\sigma$ is a common sigmoid function, ensures feasible design densities.
This feasibility property naturally extends to higher-order approximations, making the SiMPL method robust and versatile for advanced applications.
We numerically observe mesh-independent convergence of the SiMPL method, which is not entirely surprising as the method can also be derived rigorously at the infinite-dimensional function space level \cite{simplmath}.
An implementation of the SiMPL method is publicly available as an official MFEM example (Example 37) \cite{andrej2024mfem} and all of the code to reproduce our findings can be found in \cite{mfem_simpl}.
The method can be developed further to handle a larger number of design constraints than considered here. One possible extension is to utilize Lagrangian multipliers for the additional constraints. In this case, the density field can be updated by applying the SiMPL method to an augmented objective function featuring Lagrange multiplier terms that can be updated using standard augmented Lagrangian update rules.

\section*{Acknowledgments}

This work was performed under the auspices of the U.S. Department of Energy by Lawrence Livermore National Laboratory under Contract DE-AC52-07NA27344 and the LLNL-LDRD Program under Project tracking Nos.\ 22-ERD-009 and 25-ERD-030. Release number LLNL-JRNL-871320.
DK, BL, and BK were partially supported by the LLNL-LDRD Program under Project Tracking Nos.\ 22-ERD-009 and 25-ERD-030.
DK and BK were also supported in part by the U.S.\ Department of Energy Office of Science Early Career Research Program under Award Number DE-SC0024335.
\section*{Conflict of interest}
On behalf of all authors, the corresponding author states that there is no conflict of interest.
\section*{Mandatory replication of results}
All numerical experiments presented in this paper can be fully replicated using the code available at \cite{mfem_simpl}, corresponding to commit \texttt{022954b}.
\begin{appendices}
\section{Appendix}
\label{sec:appendix}
\subsection{Derivation of the gradient}

We briefly sketch the key steps to deriving a useful formula for the gradient of the objective function in the reduced setting in which we set $\widetilde{F}(\tilde{\rho}) := \hat{F}(\tilde{\rho},u(\tilde{\rho}))$. Afterwards, we replace $\tilde{\rho}$ by the mapping $\tilde{\rho}(\rho)$ and use the chain rule to derive the gradient formulae \eqref{eqs:grad} for $F(\rho) := \widetilde{F}(\tilde{\rho}(\rho))$. These are standard computations well-known in the literature on PDE-constrained optimization, e.g., \cite{Hinze2009Optimization}. We assume throughout that all functionals and operators are sufficiently smooth to allow these computations. The necessary regularity results can be rigorously derived using, e.g., elliptic regularity theory \cite{Bensoussan2002Regularity}. The canonical embeddings, e.g., that take $\rho$ into the dual space $(H^1(\Omega))^*$ in the filter equation, are left off throughout.

Given $\tilde{\rho} \in H^1(\Omega)$ and a perturbation $\delta \tilde{\rho} \in H^1(\Omega)$, we first note that the directional derivative of $\widetilde{F}$ at $\tilde{\rho}$ in direction $\delta \tilde{\rho}$ takes the form
\[
\widetilde{F}'(\tilde{\rho};\delta \tilde{\rho}) 
= 
\tilde{F}'_1(\tilde{\rho}, u(\tilde{\rho}))\delta \tilde{\rho} + 
\tilde{F}'_2(\tilde{\rho}, u(\tilde{\rho}))u'(\tilde{\rho})\delta \tilde{\rho}.
\]
Here, $\tilde{F}'_i$ for $i=1,2$ represents the partial derivative of $\tilde{F}$ with 
respect to the first and second components. 
Using the calculus of adjoints, we can make the latter term more explicit.  To start, the sensitivity equation associated with $u'(\tilde{\rho}) \delta \tilde{\rho} =: d$
is given by
\begin{equation}
  \label{eq:sen-elasticity}
  \begin{aligned}
    & \int_\Omega\Big(r(\tilde{\rho})\mathsf{C}\varepsilon(d) +  \delta \tilde{\rho} r'(\tilde{\rho})\mathsf{C}\varepsilon(u)\Big):\varepsilon(v)\dd x  = 0,
   \end{aligned}
\end{equation}
for all test functions $v\in V$. For readability, we write this in operator form:
\[
A(\tilde{\rho}) d + [A'(\tilde{\rho})\delta \tilde{\rho}]u = 0.
\]
In order words, $d = -A(\tilde{\rho})^{-1}[A'(\tilde{\rho})\delta \tilde{\rho}]u$. Note that these are merely consequences of the classical implicit function theorem. Noting that $A(\tilde{\rho})$ is a self-adjoint (symmetric) operator, we introduce the adjoint variable $\lambda := A(\tilde{\rho})^{-1} \tilde{F}'_2(\tilde{\rho}, u(\tilde{\rho}))$. We now have
\begin{equation*}
  \begin{aligned}
    \tilde{F}'_2(\tilde{\rho}, u(\tilde{\rho}))u'(\tilde{\rho})\delta \tilde{\rho}
    =\\
    -\langle \tilde{F}'_2(\tilde{\rho}, u(\tilde{\rho})),  A(\tilde{\rho})^{-1}[A'(\tilde{\rho})\delta \tilde{\rho}]u\rangle
    =\\
    -\langle A(\tilde{\rho})^{-1} \tilde{F}'_2(\tilde{\rho}, u(\tilde{\rho})), [A'(\tilde{\rho})\delta \tilde{\rho}]u\rangle
    =\\
    -\langle \lambda, [A'(\tilde{\rho})\delta \tilde{\rho}]u\rangle.
    \end{aligned}
\end{equation*}
More concretely,
\[
\widetilde{F}'(\tilde{\rho};\delta \tilde{\rho}) =  
\tilde{F}'_1(\tilde{\rho}, u(\tilde{\rho}))\delta \tilde{\rho}  - \langle \lambda, [A'(\tilde{\rho})\delta \tilde{\rho}]u\rangle,
\]
where $\lambda$ solves the adjoint equation~\eqref{eq:lambda}:
\begin{equation*}
  \int_\Omega\Big(r(\tilde{\rho})\mathsf{C}\varepsilon(\lambda)\Big):\varepsilon(v)\dd x=\langle \tilde{F}'_2(\tilde{\rho},u(\tilde{\rho})),v\rangle
\end{equation*}
for all $v\in V$ and 
\[
\langle \lambda, [A'(\tilde{\rho})\delta \tilde{\rho}]u\rangle = \int_\Omega \big(\delta \tilde{\rho}r'(\tilde{\rho})\mathsf{C}\varepsilon(u)\big):\varepsilon(\lambda)\dd x.
\]
Finally, if we denote the linear operator associated with the filter PDE~\eqref{eq:conti-filter-eq} by $L_{\epsilon}$, then using $\tilde{\rho}(\rho) = L^{-1}_{\epsilon}(\rho)$, it follows from the chain rule that
\[
F'(\rho)\delta \rho = \langle L^{-1}_{\epsilon} \widetilde{F}'(\tilde{\rho}(\rho)),\delta \rho\rangle.
\]
This is yields \eqref{eq:gtilde} and \eqref{eq:derivative}.

\end{appendices}
\bibliography{main.bbl}

\end{document}